\documentclass[twoside,11pt]{article}
\usepackage{jmlr2e}
\usepackage{array}
\usepackage{times}
\usepackage{amsmath}
\usepackage{url}
\usepackage{algorithm}
\usepackage{diagbox}
\usepackage{bbm}
\usepackage{bm}
\usepackage{diagbox}
\usepackage{graphicx}
\usepackage{subcaption}
\usepackage{epstopdf}
\usepackage{booktabs}
\usepackage{courier}
\usepackage{color}
\def\beq{\begin{eqnarray}}
\def\eeq{\end{eqnarray}}
\def\noi{\noindent}
\def\nn{\nonumber}
\def\la{\langle}
\def\ra{\rangle}

\newcommand{\bbb}[1]{\boldsymbol{\mathbf{#1}}}

\usepackage{listings}
\def\H{\mathcal{H}}

\def\J{\mathcal{J}}
\def\La{\mathcal{L}}

\def\zero{\textbf{0}}

\def\ghs{\hspace{0.3cm}} 

\newcommand{\figlabsize}{\fontsize{10}{10}\selectfont}

\usepackage{multirow}
\usepackage{rotating}

\def\imgwidlogloss{0.32\textwidth}
\def\imgheilogloss{3.5cm}
\def\imgheiseg{3.5cm}

\def\imgwidmrf{0.32\textwidth}
\def\imgheimrf{3.4cm}

\def\imgwidl0tv{0.32\textwidth}
\def\imgheil0tv{3.3cm}
\jmlrheading{1}{2016}{1-48}{4/00}{10/00}{Ganzhao Yuan and Bernard Ghanem}

\ShortHeadings{Sparsity Constrained Minimization}{Yuan and Ghanem}
\firstpageno{1}
\let\cite\citep
\definecolor{colortwo}{RGB}{128,128,127}
\definecolor{colorone}{RGB}{218,218,217}

\begin{document}

\title{Sparsity Constrained Minimization via Mathematical Programming with Equilibrium Constraints}
\author{\name Ganzhao Yuan \email yuanganzhao@gmail.com \\
       \addr School of Data \& Computer Science, Sun Yat-sen University (SYSU), P.R. China
       \AND
       \name Bernard Ghanem \email bernard.ghanem@kaust.edu.sa \\
       \addr King Abdullah University of Science \& Technology (KAUST), Saudi Arabia}
\editor{}

\maketitle

\begin{abstract}

Sparsity constrained minimization captures a wide spectrum of applications in both machine learning and signal processing. This class of problems is difficult to solve since it is NP-hard and existing solutions are primarily based on Iterative Hard Thresholding (IHT). In this paper, we consider a class of continuous optimization techniques based on Mathematical Programs with Equilibrium Constraints (MPECs) to solve general sparsity constrained problems. Specifically, we reformulate the problem as an equivalent biconvex MPEC, which we can solve using an exact penalty method or an alternating direction method. We elaborate on the merits of both proposed methods and analyze their convergence properties. Finally, we demonstrate the effectiveness and versatility of our methods on several important problems, including feature selection, segmented regression, trend filtering, MRF optimization and impulse noise removal. Extensive experiments show that our MPEC-based methods outperform state-of-the-art techniques, especially those based on IHT. 

\noi \textbf{Keywords:} Sparsity Constrained Optimization, Exact Penalty Method, Alternating Direction Method, MPEC, Convergence Analysis
\end{abstract}

\section{Introduction}
In this paper, we mainly focus on the following generalized sparsity constrained minimization problem:
\beq \label{eq:l0}
\min_{\bbb{x}}~ f(\bbb{x}),~s.t.~\|\bbb{Ax}\|_0 \leq k
\eeq
\noi where $\bbb{x}\in \mathbb{R}^{n}$, $\bbb{A}\in \mathbb{R}^{m\times n}$ and $k~ (0<k<m)$ is a positive integer. $\|\cdot\|_0$ is a function that counts the number of nonzero elements in a vector. To guarantee convergence, we assume that $f(\cdot)$ is convex and $L$-Lipschitz continuous (but not necessarily smooth), $\bbb{A}$ has right inverse (i.e. $rank(\bbb{A})=m$), and there always exists bounded solutions to Eq (\ref{eq:l0}).




The optimization in Eq (\ref{eq:l0}) describes many applications of interest in both machine learning and signal processing, including compressive sensing \cite{donoho2006compressed}, Dantzig selector \cite{candes2007dantzig}, feature selection \cite{ng2004feature}, trend filtering \cite{KimKBG09}, image restoration \cite{YuanG15,Dong2013}, sparse classification and boosting \cite{WestonEST03,XiangXHR09}, subspace clustering \cite{elhamifar2009sparse}, sparse coding \cite{LeeBRN06,mairal2010online,BaoJQS16}, portfolio selection, image smoothing \cite{XuLXJ11}, sparse inverse covariance estimation \cite{friedman2008sparse,yuan2010high}, blind source separation\cite{li2004analysis}, permutation problems \cite{FogelJBd15,jiang2015lp}, joint power and admission control \cite{LiuDM15}, Potts statistical functional \cite{weinmann20151}, to name a few. Moreover, we notice that many binary optimization problems \cite{WangSH13} can be reformulated as an $\ell_0$ norm optimization problem, since $\bbb{x}\in\{-1,+1\}^n \Leftrightarrow \|\bbb{x}-\bbb{1}\|_0 + \|\bbb{x}+\bbb{1}\|_0\leq n$. In addition, it can be rewritten in a more compact form as $\|\bbb{Ax}-\bbb{b}\|_0\leq n$ with $\bbb{A}=[\bbb{I}_n, \bbb{I}_n]^T \in \mathbb{R}^{2n \times n}$ and $\bbb{b}=[-\bbb{1}^T, \bbb{1}^T]^T\in \mathbb{R}^{2n}$.

A popular method to solve Eq (\ref{eq:l0}) is to consider its $\ell_1$ norm convex relaxation, which simply replaces the $\ell_0$ norm function\footnote{Strictly speaking, $\|\cdot\|_0$ is not a norm (but a pseudo-norm) since it is not positive homogeneous.} with the $\ell_1$ norm function. When $\bbb{A}$ is identity and $f(\cdot)$ is a least squares objective function, it has been proven that when the sampling matrix in $f(\cdot)$ satisfies certain incoherence conditions and $\bbb{x}$ is sparse at its optimal solution, the problem can be solved exactly by this convex method \cite{CandesT05}. However, when these two assumptions do not hold, this convex method can be unsatisfactory \cite{LiuDM15,YuanG15}.

Recently, another breakthrough in sparse optimization is the multi-stage convex relaxation method \cite{zhang2010analysis,candes2008enhancing}. In the work of \cite{zhang2010analysis}, the local solution obtained by this method is shown to be superior to the global solution of the standard $\ell_1$ convex method, which is used as its initialization. However, this type of method seeks an approximate solution to the sparse \emph{regularized} problem with a smooth objective, and cannot solve the sparsity \emph{constrained} problem \footnote{Some may solve the constrained problem by tuning the continuous parameter of regularized problem, but it can be hard when $m$ is large.} considered here. Our method uses the same convex initialization strategy, but it is applicable to general $\ell_0$ norm minimization problems.

\vspace{3pt}\noi \textbf{Challenges and Contributions:} We recognize three main challenges hindering existing work. \textbf{(a)} The general problem in Eq (\ref{eq:l0}) is NP-hard. There is little hope of finding the global minimum efficiently in all instances. In order to deal with this issue, we reformulate the $\ell_0$ norm minimization problem as an equivalent augmented optimization problem with a bilinear equality constraint using a variational characterization of the $\ell_0$ norm function. Then, we propose two penalization/regularization methods (exact penalty and alternating direction) to solve it. The resulting algorithms seek a desirable exact solution to the original optimization problem. \textbf{(b)} Many existing convergence results for non-convex $\ell_0$ minimization problems tend to be limited to unconstrained problems or inapplicable to constrained optimization. We carefully analyze both our proposed algorithms and prove that they always converge to a first-order KKT point. To the best of our knowledge, this is the first attempt to solve general $\ell_0$ constrained minimization with guaranteed convergence. \textbf{(c)} Exact methods (e.g. methods using IHT) can produce unsatisfactory results, while approximation methods such as Schatten's $\ell_p$ norm method \cite{GeJY11} and re-weighted $\ell_1$ \cite{candes2008enhancing} fail to control the sparsity of the solution. In comparison, experimental results show that our MPEC-based methods outperform state-of-the-art techniques, especially those based on IHT. This is consistent with our new technical report \cite{yuan2016binary} which shows that IHT method often presents sub-optimal performance in binary optimization.





%
\vspace{3pt}\noi \textbf{Organization and Notations:} This paper is organized as follows. Section \ref{sect:related} provides a brief description of the related work. Section \ref{sect:epm} and Section \ref{sect:adm} present our proposed exact penalty method and alternating direction method optimization framework, respectively. Section \ref{sect:disc} discusses some features of the proposed two algorithms. Section \ref{sect:exp} summarizes the experimental results. Finally, Section \ref{sect:conc} concludes this paper. Throughout this paper, we use lowercase and uppercase boldfaced letters to denote real vectors and matrices respectively. We use $\la \bbb{x},\bbb{y}\ra$ and $\bbb{x} \odot \bbb{y}$  to denote the Euclidean inner product and elementwise product between $\bbb{x}$ and $\bbb{y}$. ``$\triangleq$'' means define. $\sigma(\bbb{A})$ denotes the smallest singular value of $\bbb{A}$. 

\begin{table*}[!t]
\begin{center}
\caption{$\ell_0$ norm optimization techniques.}\label{tab:l0opt}
\scalebox{0.825}{\begin{tabular}{|p{8.8cm}|p{8.7cm}|}
  \hline
  Method and Reference & Description       \\
  \hline
  greedy descent methods \text{\cite{mallat1993matching}}  & \text{ only for smooth (typically quadratic) objective}   \\
  \hline
    $\ell_1$ norm relaxation \text{\cite{CandesT05}}  & \text{$\|\bbb{x}\|_0 \approx \|\bbb{x}\|_1$}  \\
  \hline
   \text{ $k$-support norm relaxation \cite{ArgyriouFS12}}  & \text{$\|\bbb{x}\|_0 \approx \|\bbb{x}\|_{\text{k-sup}} \triangleq \max_{\bbb{0}< \bbb{v}\leq \bbb{1},~\la \bbb{v} ,\bbb{1}\ra \leq k}~(\sum_{i}{\bbb{x}_i^2}/{\bbb{v}_i})^{1/2}$}  \\
  \hline
   \text{ $k$-largest norm relaxation \cite{yu2011adversarial}}  & \text{$\|\bbb{x}\|_0 \approx \|\bbb{x}\|_{\text{k-lar}} \triangleq \max_{-\bbb{1}\leq \bbb{v}\leq \bbb{1},~\|\bbb{v}\|_{1}\leq k}~\la \bbb{v}, \bbb{x} \ra $}  \\
  \hline
  \text{ SOCP convex relaxation \cite{Chan2007} } & \text{$\|\bbb{x}\|_0\leq k \Rightarrow \|\bbb{x}\|_1 \leq \sqrt{k}\|\bbb{x}\|_2$}   \\
  \hline
   \text{ SDP convex relaxation \cite{Chan2007} } & \text{$\|\bbb{x}\|_0\leq k  \Rightarrow \bbb{X} = \bbb{xx}^T, ~\|\bbb{X}\|_1\leq k tr(\bbb{X})$} \\ 
  \hline
   \text{Schatten $\ell_p$ approximation \cite{GeJY11} } & \text{$\|\bbb{x}\|_0 \approx \|\bbb{x}\|_p$}  \\
  \hline
   \text{ re-weighted $\ell_1$ approximation \cite{candes2008enhancing}} & \text{$\|\bbb{x}\|_0 \approx \la \textbf{1}, \ln(|\bbb{x}|+\epsilon)\ra $} \\
  \hline
   \text{ $\ell_{\text{1-2}}$ DC approximation \cite{YinLHX15} } & \text{$\|\bbb{x}\|_0 \approx \|\bbb{x}\|_1 - \|\bbb{x}\|_2 $} \\
  \hline
  \text{ $0\text{-}1$ mixed integer programming \cite{bienstock1996computational} }& \text{$\{\|\bbb{x}\|_0:\|\bbb{x}\|_{\infty} \leq \lambda \} \Leftrightarrow \{\min_{\bbb{v}\in\{0,1\}}~\la \bbb{1},\bbb{v} \ra: |\bbb{x}|\leq \lambda \bbb{v}\}$} \\
  \hline
 \text{ iterative hard shreadholding \cite{beck2013sparsity} } & \text{$0.5\|\bbb{x}-\bbb{x}'\|_2^2,~s.t.~\|\bbb{x}\|_0\leq k$}  \\
  \hline
  \text{ non-separable MPEC \hspace*{\fill}[This paper]} &\text{$\|\bbb{x}\|_0 = \min_{\bbb{u}}~\|\bbb{u}\|_1, ~s.t. ~\|\bbb{x}\|_1 =  \la \bbb{x},\bbb{u} \ra,~-\bbb{1}\leq \bbb{u} \leq \bbb{1}$} \\
  \hline
  \text{ separable MPEC \cite{YuanG15} \hspace*{\fill}[This paper]} & \text{$\|\bbb{x}\|_0 = \min_{\bbb{v}}~\la \bbb{1}, \bbb{1}-\bbb{v}\ra, ~s.t.~|\bbb{x}|\odot \bbb{v}=\bbb{0},~\bbb{0}\leq \bbb{v}\leq \bbb{1}$} \\
  \hline
\end{tabular}}
\end{center}
\end{table*}

\section{Related Work} \label{sect:related}
There are mainly four classes of $\ell_0$ norm minimization algorithms in the literature: \textbf{(i)} greedy descent methods, \textbf{(ii)} convex approximate methods, \textbf{(iii)} non-convex approximate methods, and \textbf{(iv)} exact methods. We summarize the main existing algorithms in Table \ref{tab:l0opt}.

\textbf{Greedy descent methods.} They have a monotonically decreasing property and optimality guarantees in some cases \cite{tropp2004greed}, but they are limited to solving problems with smooth objective functions (typically the square function). \textbf{(a)} Matching pursuit \cite{mallat1993matching} selects at each step one atom of the variable that is the most correlated with the residual. \textbf{(b)} Orthogonal matching pursuit \cite{tropp2007signal} uses a similar strategy but also creates an orthonormal set of atoms to ensure that selected components are not introduced in subsequent steps. 
\textbf{(c)} Gradient pursuit \cite{blumensath2008gradient} is similar to matching pursuit, but it updates the sparse solution vector at each iteration with a directional update computed based on gradient information. \textbf{(d)} Gradient support pursuit \cite{BahmaniRB13} iteratively chooses the index with the largest magnitude as the pursuit direction while maintaining the stable restricted Hessian property. \textbf{(e)} Other greedy methods has been proposed, including basis pursuit \cite{chen1998atomic}, regularized orthogonal matching pursuit \cite{needell2009uniform}, compressive sampling matching pursuit \cite{needell2009cosamp}, and forward-backward greedy method \cite{zhang2011adaptive,RaoSW15}.


\textbf{Convex approximate methods.} They seek convex approximate reformulations of the $\ell_0$ norm function. \textbf{(a)} The $\ell_1$ norm convex relaxation\cite{CandesT05} provides a convex lower bound of the $\ell_0$ norm function in the unit $\ell_{\infty}$ norm. It has been proven that under certain incoherence assumptions, this method leads to a near optimal sparse solution. However, such assumptions may be violated in real applications. \textbf{(b)} $k$-support norm provides the tightest convex relaxation of sparsity combined with an $\ell_2$ penalty. Moreover, it is tighter than the elastic net \cite{zou2005regularization} by exactly a factor of $\sqrt{2}$ \cite{ArgyriouFS12,mcdonald2014spectral}. \textbf{(c)} $k$-largest norm provides the tightest convex Boolean relaxation in the sense of minimax game \cite{yu2011adversarial,PilanciWG15,Dattorro2011}. \textbf{(d)} Other convex Second-Order Cone Programming (SOCP) relaxations and Semi-Definite Programming (SDP) relaxations \cite{Chan2007} have been considered in the literature.

\textbf{Non-convex approximate methods.} They seek non-convex approximate reformulations of the $\ell_0$ norm function. \textbf{(a)} Schatten $\ell_p$ norm with $p \in (0,1)$ was considered by \cite{GeJY11} to approximate the discrete $\ell_0$ norm function. It results in a local gradient Lipschitz continuous function, to which some smooth optimization algorithms can be applied. \textbf{(b)} Re-weighted $\ell_1$ norm \cite{candes2008enhancing,zhang2010analysis,zou2006adaptive} minimizes the first-order Taylor series expansion of the objective function iteratively to find a local minimum. It is expected to refine the $\ell_0$ regularized problem, since its first iteration is equivalent to solving the $\ell_1$ norm problem. \textbf{(c)} $\ell_{\text{1-2}}$ DC (difference of convex) approximation \cite{YinLHX15} was considered for sparse recovery. Exact stable sparse recovery error bounds were established under a restricted isometry property condition. \textbf{(d)} Other non-convex surrogate functions of the $\ell_0$ function have been proposed, including Smoothly Clipped Absolute Deviation (SCAD) \cite{fan2001variable}, Logarithm \cite{friedman2012fast}, Minimax Concave Penalty (MCP) \cite{zhang2010nearly}, Capped $\ell_1$ \cite{zhang2010analysis}, Exponential-Type \cite{Gao2011}, Half-Quadratic \cite{geman1995nonlinear}, Laplace \cite{trzasko2009highly}, and MC+ \cite{mazumder2012sparsenet}. Please refer to \cite{lu2014generalized,GongZLHY13} for a detailed summary and discussion.


\textbf{Exact methods.} Despite the success of approximate methods, they are unappealing in cases when an exact solution is required. Therefore, many researchers have sought out exact reformulations of the $\ell_0$ norm function. \textbf{(a)} $0\text{-}1$ mixed integer programming \cite{bienstock1996computational,bertsimas2015best} assumes the solution has bound constraints. It can be solved by a tailored branch-and-bound algorithm, where the cardinality constraint is replaced by a surrogate constraint. \textbf{(b)} Hard thresholding \cite{LuZ13,beck2013sparsity,yuan2014gradient,jain2014iterative} iteratively sets the smallest (in magnitude) values to zero in a gradient descent format. It has been incorporated into the Quadratic Penalty Method (QPM) \cite{LuZ13} and Mean Doubly Alternating Direction Method \cite{Dong2013} (MD-ADM). However, we found they often converge to unsatisfactory results in practise. \textbf{(c)} The MPEC reformulation \cite{YuanG15} considers the complementary system of the $\ell_0$ norm problem by introducing additional dual variables and minimizing the complimentary error of the MPEC problem. However, it is limited to $\ell_0$ norm regularized problems and the convergence results are weak.


From above, we observe that existing methods either produce approximate solutions or are limited to smooth objectives $f(.)$ and when $\mathbf{A}=\mathbf{I}$. The only existing general purpose exact method for solving Eq (\ref{eq:l0}) is the quadratic penalty method \cite{LuZ13} and the mean doubly alternating direction method \cite{Dong2013}. However, they often produce unsatisfactory results in practise. 
The unappealing shortcomings of existing solutions and renewed interests in MPECs \cite{YuanG15,feng2013complementarity} motivate us to design new MPEC-based algorithms and convergence results for sparsity constrained minimization.

\section{Exact Penalty Method} \label{sect:epm}

In this section, we present an exact penalty method \cite{luo1996mathematical,hu2004convergence,kadrani2009new,BiLP14} for solving the problem in Eq (\ref{eq:l0}). This method is based on an equivalent non-separable MPEC reformulation of the $\ell_0$ norm function.

\subsection{ Non-Separable MPEC Reformulation}
First of all, we present a new non-separable MPEC reformulation.

 \begin{lemma}\label{lemma:non:sep:mpec}
 For any given $\bbb{x}\in \mathbb{R}^{n}$, it holds that
  \beq\label{eq:non:sep:mpec}
  \|\bbb{x}\|_0=\min_{\bbb{-1}\le \bbb{u}\le \bbb{1}}~\|\bbb{u}\|_1,~s.t.~\|\bbb{x}\|_1 = \la \bbb{u}, \bbb{x} \ra,
 \eeq
 and $\bbb{u}^*={\rm sign}(\bbb{x})$ is the unique optimal solution of Eq(\ref{eq:non:sep:mpec}). Here, the standard signum function sign is applied componentwise, and ${\rm sign}(0)=0$.

\begin{proof}

It is not hard to validate that the $\ell_0$ norm function of $\bbb{x}\in \mathbb{R}^m$ can be repressed as the following minimization problem over $\bbb{u}\in \mathbb{R}^m$:
  \beq \label{eq:mpec2:proof}
\|\bbb{x}\|_0  = \min_{-\bbb{1} \leq \bbb{u} \leq \bbb{1}}~ \|\bbb{u}\|_1,~s.t.~|\bbb{x}| = \bbb{u} \odot \bbb{x}
  \eeq
\noi Note that the minimization problem in Eq (\ref{eq:mpec2:proof}) can be decomposed into $m$ subproblems. When $\bbb{x}_i=0$, the optimal solution in position $i$ will be achieved at $\bbb{u}_i=0$ by minimization; when $\bbb{x}_i>0$ ($\bbb{x}_i<0$, respectively), the optimal solution in position $i$ will be achieved at $\bbb{u}_i=1$ ($\bbb{u}_i=-1$, respectively) by constraint. In other words, $\bbb{u}^*={\rm sign}(\bbb{x})$ will be achieved for Eq (\ref{eq:mpec2:proof}), leading to $\|\bbb{x}\|_0=\|{\rm sign}(\bbb{x})\|_1$.

We now focus on Eq (\ref{eq:mpec2:proof}). Due to the box constraints $-\bbb{1} \leq \bbb{u} \leq \bbb{1}$, it always holds that $|\bbb{x}| - \bbb{u}\odot \bbb{x} \geq \bbb{0}$. Therefore, the error generated by the difference of $|\bbb{x}|$ and $\bbb{u}\odot \bbb{x}$ is summarizable. We naturally have the following results:
\beq
|\bbb{x}| - \bbb{u} \odot \bbb{x} = \bbb{0}~ &\Leftrightarrow& ~\la \bbb{1}, |\bbb{x}| - \bbb{u}\odot \bbb{x} \ra =0 \Leftrightarrow\|\bbb{x}\|_1 = \la \bbb{u},\bbb{x} \ra \nn
\eeq
\end{proof}
\end{lemma}

We remark that similar conclusions of this lemma have been appeared in \cite{BiLP14,Bi2014,feng2013complementarity,BiP17} and we present a different reformulation here.

Using Lemma \ref{lemma:non:sep:mpec}, we can rewrite Eq (\ref{eq:l0}) in an equivalent form as follows.
\beq \label{eq:l0:mpec:nonsep}
\min_{\bbb{x},\bbb{u}}~f(\bbb{x})+I(\bbb{u}),~s.t.~\|\bbb{Ax}\|_1 = \la \bbb{Ax} , \bbb{u} \ra
\eeq
\noi where $I(\bbb{u})$ is an indicator function on $\Theta$ with
\beq
I(\bbb{u}) {\scriptsize = \begin{cases}
0,~\bbb{u} \in \Theta\\
\infty,~\bbb{u} \notin \Theta,\\
\end{cases}}\Theta \triangleq \{\bbb{u}|~\|\bbb{u}\|_1\leq   k, ~\bbb{-1}\leq \bbb{u}\leq \bbb{1}\}. \nn
\eeq
\noi


\begin{algorithm} [!h]
\caption{\label{alg:epm} {MPEC-EPM: An Exact Penalty Method for Solving MPEC Problem (\ref{eq:l0:mpec:nonsep})}}
\raggedright (S.0) Initialize $\bbb{x}^0= 0 \in \mathbb{R}^{n}$, $\bbb{u}^0 = \bbb{0} \in \mathbb{R}^{m}$, $\rho^{0} >0$. Set $t=0$ and $\mu>0$.

\raggedright (S.1) Solve the following $\bbb{x}$-subproblem:
     \beq \label{eq:subprob:epm:x}
     \bbb{x}^{t+1}= \mathop{\arg\min}_{\bbb{x}}~\J(\bbb{x},\bbb{u}^{t}) + \frac{\mu}{2} \|\bbb{x}-\bbb{x}^t\|_2^2
     \eeq

\raggedright (S.2) Solve the following $\bbb{u}$-subproblem:
     \beq \label{eq:subprob:epm:u}
      \bbb{u}^{t+1} =\mathop{\arg\min}_{ \bbb{u}}~\J(\bbb{x}^{t+1},\bbb{u}) + \frac{\mu}{2} \|\bbb{u}-\bbb{u}^t\|_2^2
     \eeq

\raggedright (S.3) Update the penalty in every $T$ iterations:
                \beq \label{eq:update:epm:rho}
                \rho^{t+1}= \min(L/\sigma(\bbb{A}),2 \rho^{t})
                \eeq
\raggedright (S.4) Set $t:=t+1$ and then go to Step (S.1)
 \end{algorithm}

\subsection{Proposed Optimization Framework}
We now give a detailed description of our solution algorithm to the optimization in Eq (\ref{eq:l0:mpec:nonsep}). Our solution is based on the exact penalty method, which penalizes the complementary error directly by a penalty function. The resulting objective $\mathcal{J}: \mathbb{R}^{n}\times\mathbb{R}^{m}\to\mathbb{R}$ is defined in Eq (\ref{eq:main:mpec:penalty}), where $\rho$ is the tradeoff penalty parameter that is iteratively increased to enforce the constraints.
\beq \label{eq:main:mpec:penalty}
\begin{split}
\mathcal{J}_{\rho}(\bbb{x},\bbb{u}) \triangleq f(\bbb{x}) + I(\bbb{u}) + \rho ( \|\bbb{Ax}\|_1 - \la \bbb{Ax}, \bbb{u} \ra)
\end{split}
\eeq
\noi  In each iteration, $\rho$ is fixed and we use a proximal point method \cite{Attouch2010,Bolte2014} to minimize over $\bbb{x}$ and $\bbb{u}$ in an alternating fashion. Details of this exact penalty method are in Algorithm \ref{alg:epm}. Note that the parameter $T$ is the number of inner iterations for solving the bi-convex problem. We make the following observations about the algorithm.

\vspace{1pt}\noi \bbb{(a) Initialization}. We initialize $\bbb{u}^0$ to $\bbb{0}$. This is for the sake of finding a reasonable local minimum in the first iteration, as it reduces to a convex $\ell_1$ norm minimization problem for the $\bbb{x}$-subproblem. 


\vspace{1pt}\noi \bbb{(b) Exact property}. The key feature of this method is the boundedness of the penalty parameter $\rho$ (see Theorem \ref{theorem:1}). Therefore, we terminate the optimization when the threshold is reached (see Eq (\ref{eq:update:epm:rho})). This distinguishes it from the quadratic penalty method \cite{LuZ13}, where the penalty may become arbitrarily large for non-convex problems.

\vspace{1pt}\noi \bbb{(c) $\bbb{u}$-Subproblem}. Variable $\bbb{u}$ in Eq (\ref{eq:subprob:epm:u}) is updated by solving the following problem:
\beq \label{eq:u:subp}
\begin{split}
\bbb{u}^{t+1} =\mathop{\arg\min}_{\bbb{-1}\leq \bbb{u}\leq \bbb{1}}   ~ \frac{1}{2} \| \bbb{u} - \bbb{z}\|_2^2,~s.t.~\|\bbb{u}\|_1\leq k,
\end{split}
\eeq
\noi where $\bbb{z}=\rho \bbb{Ax}^{t+1}/\mu$ and $\mu$ is a non-negative proximal constant. 
Due to the symmetry in the objective and constraints, we can without loss of generality assume that $\bbb{z}\geq 0$, and flip signs of the resulting solution at the end. As such, Eq (\ref{eq:u:subp}) can be solved by
\beq \label{eq:u:subp2}
\begin{split}
\bbb{\bar{u}}^{t+1} =\mathop{\arg\min}_{\bbb{0}\leq \bbb{u}\leq \bbb{1}}   ~ \frac{1}{2} \| \bbb{u} - |\bbb{z}|\|_2^2,~s.t.~ \la \bbb{u}, \bbb{1} \ra \leq k,
\end{split}
\eeq
\noi where the solution in Eq (\ref{eq:u:subp}) can be recovered via $\bbb{u}^{t+1}=sign(\bbb{z}) \odot \bbb{\bar{u}}^{t+1}$. Eq (\ref{eq:u:subp2}) can be solved \emph{exactly} in $n\log(n)$ time using the breakpoint search algorithm \cite{HelgasonKL80}. Note that this algorithm includes the simplex projection \cite{DuchiSSC08} as a special case. For completeness, we also include a Matlab code in \textbf{Appendix} \ref{supp:breakpoint}.


\noi \bbb{(d) $\bbb{x}$-Subproblem}. Variable $\bbb{x}$ in Eq (\ref{eq:subprob:epm:x}) is updated by solving a $\ell_1$ norm convex problem which has no closed-form solution. However, it can be solved using existing $\ell_1$ norm optimization solvers such as classical or linearized ADM \cite{HeY12} with convergence rate $\mathcal{O}(1/t)$.




\subsection{Theoretical Analysis}
We present the convergence analysis of the exact penalty method. Our main novelties are establishing the inversely proportional relationship between the weight of $\ell_1$ norm and the sparsity of $\bbb{Ax}$ (in the proof of exactness) and qualifying the sparsity upper bound for $\la \bbb{w}, |\bbb{Ax}|\ra$ (in the proof of convergence rate). Our results make use of the specified structure of the $\ell_1$ norm problem.

The following lemma is useful in our convergence analysis.
\begin{lemma}\label{lemma:hhhh}
Assume that $\bbb{A}\in \mathbb{R}^{m\times n}$ has right inverse, $h(\cdot)$ is convex and $L$-Lipschitz continuous, and $\bbb{w}$ is a non-negative vector. We have the following results: (i) It always holds that $\|\bbb{A}^T\bbb{c}\|_2 \geq \sigma(\bbb{A})\|\bbb{c}\|_{2}$ for all $\bbb{c}\in \mathbb{R}^m$, where $\sigma(\bbb{A})$ is the smallest singular value of $\bbb{A}$. (ii) The optimal solution of the following optimization problem:
\beq \label{eq:sub:Ax}
\bbb{x}^* = \arg \min_{\bbb{x}} ~h(\bbb{x}) + \la \bbb{w}, |\bbb{Ax}|\ra
\eeq
\noi will be achieved with $|\bbb{Ax}^*|_i=0$ when $\bbb{w}_i > L/ \sigma(\bbb{A})$ for all $i\in[m]$. (iii) Moreover, if $\|\bbb{x}^*\|\leq \delta$, it always holds that $\la \bbb{w}, |\bbb{Ax}^*|\ra \leq \delta L $.

\begin{proof}

(i) We now prove the first part of this lemma. (i) We denote $\bbb{Z}=\bbb{A}\bbb{A}^T - (\sigma(\bbb{A}))^2 \bbb{I} \in \mathbb{R}^{m\times m}$. Since $\bbb{A}$ has right inverse, we have $\sigma(\bbb{A})>0$ and $\bbb{Z}\succeq \bbb{0}$. We let $\bbb{Z}=\bbb{L}\bbb{L}^T$ and have the following inequalities:
\beq
\|\bbb{A}^T\bbb{c}\|_2^2 &=& \la \bbb{A}\bbb{A}^T ,\bbb{cc}^T \ra\nn\\
&=& \la \bbb{Z} + (\sigma(\bbb{A}))^2 \bbb{I} ,\bbb{cc}^T\ra \nn \\
&=& \|\bbb{L}^T\bbb{c}\|_2^2 + \la  (\sigma(\bbb{A}))^2 \bbb{I}, \bbb{cc}^T\ra\nn\\
&\geq& 0 + \la  (\sigma(\bbb{A}))^2 \|\bbb{c}\|^2 \nn
\eeq
\noi Taking the square root of both sides, we have
\beq \label{eq:jjjjj}
\|\bbb{A}^T\bbb{c}\|_2 \geq \sigma(\bbb{A})\|\bbb{c}\|_{2}.
\eeq

(ii) We now prove the second part of this lemma. The convex optimization problem in Eq (\ref{eq:sub:Ax}) is equivalent to the following minimax saddle point problem:
\beq\label{eq:sub:saddle}
(\bbb{x}^*,\bbb{y}^*)=\arg\min_{\bbb{x}}~\arg\max_{-\bbb{1}\leq \bbb{y}\leq \bbb{1}}~h(\bbb{x}) + \la \bbb{y}\odot \bbb{w}, \bbb{Ax}\ra\nn
\eeq
\noi According to the optimality with respect to $\bbb{y}$, we have the following result:
\beq \label{eq:optimal:y}
\bbb{y}^*_i=1  &\Rightarrow& (\bbb{Ax}^*)_i>0,\nn\\
|\bbb{y}^*_i|<1 &\Rightarrow& (\bbb{Ax}^*)_i=0,\\
\bbb{y}^*_i=-1  &\Rightarrow& (\bbb{Ax}^*)_i<0\nn
\eeq
\noi According to the optimality with respect to $\bbb{x}$, we have the following result:
\beq \label{eq:optimal:x}
h'(\bbb{x}) + \bbb{A}^T(\bbb{y}\odot \bbb{w})=\bbb{0}
\eeq
\noi where $h'(\bbb{x})$ denotes the subgradient of $h(\cdot)$ with respect to $\bbb{x}$.


We now prove that $|\bbb{y}_i|$ is strictly less than 1 when $\bbb{w}_i>L/\sigma(A)$ for all $i$, then we obtain $|\bbb{Ax^*}|_i=0$ due to the optimality condition in Eq (\ref{eq:optimal:y}). Our proof is as follows. We assume that $\bbb{y}_i \neq 0$, since otherwise $|\bbb{y}_i|<0$, the conclusion holds. Then we derive the following inequalities:
\beq \label{eq:y:bound1}
|\bbb{y}_i| &<& |\bbb{y}_i| \cdot \frac{  \bbb{w}_i }{L} \cdot \sigma(\bbb{A}) \nn\\
&\leq& |\bbb{y}_i| \cdot \frac{\bbb{w}_i }{L} \cdot \frac{\|\bbb{A}^T (\bbb{y}\odot \bbb{w})\|}{ \|\bbb{y}\odot \bbb{w}\| }\nn\\
&\leq& |\bbb{y}_i| \cdot \frac{\bbb{w}_i }{L} \cdot \frac{L}{ \|\bbb{y}\odot \bbb{w}\| }\nn\\
&\leq&   \frac{|\bbb{y}_i| \cdot \bbb{w}_i }{ \|\bbb{y}\odot \bbb{w}\|_{\infty} }\nn\\
&\leq&   1
\eeq
\noi where the first step uses the choice that $\bbb{w}_i>L/\sigma(A)$ and the assumption that $\bbb{y}_i \neq 0$; the second step uses (\ref{eq:jjjjj}); the third step uses the fact that $\|\bbb{A}^T(\bbb{y}\odot \bbb{w})\|=\|h'(\bbb{x})\|\leq L$ which can be derived from Eq (\ref{eq:optimal:x}) and the $L$-Lipschitz continuity of $h(\cdot)$, the fourth step uses the norm inequality that $\|\bbb{c}\|_{\infty}\leq\|\bbb{c}\|_2~\forall \bbb{c}$, the last step uses the nonnegativity of $\bbb{w}$.

We remark that such results are natural and commonplace, since it is well-known that $\ell_1$ norm induces sparsity.

(iii) We now prove the third part of this lemma. Recall that if $\bbb{x}^*$ solves Eq (\ref{eq:sub:Ax}), then it must satisfy the following variational inequality \cite{HeY12}:
\beq \label{eq:variational}
\la \bbb{w}, | \bbb{Ax}| - | \bbb{Ax}^*| \ra + \la  \bbb{x} - \bbb{x}^*, h'(\bbb{x}^*)\ra \geq 0,~\forall \bbb{x}
\eeq
\noi Letting $\bbb{x}=\bbb{0}$ in Eq (\ref{eq:variational}), using the nonnegativity of $\bbb{w}$, Cauchy-Schwarz inequality, and the Lipschitz continuity of $h(\cdot)$, we have:
\beq
\la \bbb{w},  | \bbb{Ax}^*| \ra  &\leq&  - \la  \bbb{x}^*, h'(\bbb{x}^*)\ra \nn\\
&\leq& \|\bbb{x}^*\| \cdot \|h'(\bbb{x}^*)\|\nn\\
&\leq&  \delta  L\nn
\eeq
\end{proof}
\end{lemma}

The following lemma shows that the biconvex minimization problem will lead to basic feasibility for the complementarity constraint when the penalty parameter $\rho$ is larger than a threshold.


\begin{lemma} \label{lemma:v:01}
Any local optimal solution of the minimization problem: $\min_{\bbb{x},\bbb{u}}~\J_{\rho}(\bbb{x},\bbb{u})$ will be achieved with $\|\bbb{Ax}\|_1 - \la \bbb{Ax} , \bbb{u}\ra =0$, when $\rho > L/\sigma(\bbb{A})$.

\begin{proof}

We let $\{\bbb{x},\bbb{u}\}$ be any local optimal solution of the biconvex minimization problem $\min_{\bbb{x},\bbb{u}}~\J_{\rho}(\bbb{x},\bbb{u})$. We denote $O$ as the index of the largest-k value of $|\bbb{Ax}|$, $G \triangleq \{i|(\bbb{Ax})_i>0\}$, and $S \triangleq\{j|(\bbb{Ax})_j<0\}$. Moreover, we define $I \triangleq O \cap G$, $J \triangleq O \cap S$ and $K \triangleq \{1,2,...,m\} \setminus \{I \cup J \} $.

(i) First of all, we consider the minimization problem of the penalty function with respect to $\bbb{u}$ (i.e. $\min_{\bbb{u}}~\J_{\rho}(\bbb{x},\bbb{u})$). It reduces to the following minimization problem:
\beq \label{eq:V:subproblem:1}
\min_{-\bbb{1}\leq \bbb{u} \leq \bbb{1}}~  \la -\bbb{Ax}, \bbb{u}\ra,~s.t.~\|\bbb{u}\|_1\leq k
\eeq
\noi The objective of Eq (\ref{eq:V:subproblem:1}) essentially computes the $k$-largest norm function of $\bbb{Ax}$, see \cite{WuDST14,Dattorro2011}. It is not hard to validate that the optimal solution of Eq (\ref{eq:V:subproblem:1}) can be computed as follows:
\beq \label{eq:opt:u}
\bbb{u}_i \scriptsize = \begin{cases}
+1, & i \in I\\
-1, & i\in J \\
0, & u \in K \\
\end{cases}\nn
\eeq

(ii) We now consider the minimization problem of the penalty function with respect to $\bbb{x}$ (i.e. $\min_{\bbb{x}}~\J_{\rho}(\bbb{x},\bbb{u})$). Clearly, we have that: $\la (\bbb{Ax})_I, \bbb{u}_I \ra = \| (\bbb{Ax})_I\|_1$ and $\la (\bbb{Ax})_J, \bbb{u}_J \ra = \| (\bbb{Ax})_J\|_1$. It reduces to the following optimization problem for the $\bbb{x}$-subproblem:
\beq
\min_{\bbb{x}}~f(\bbb{x})  + \rho ( \|(\bbb{Ax})_K\|_1 - \la (\bbb{Ax})_K, \bbb{u}_K \ra ).\nn
\eeq
\noi Since $\bbb{u}_K=0$, it can be further simplified as:
\beq \label{eq:KK}
\min_{\bbb{x}}~f(\bbb{x})  + \la \bbb{w},|\bbb{Ax}| \ra,~\text{where}~\bbb{w}_i=
\begin{cases}
\rho,&i\in K\\
0&i\notin K
\end{cases}
\eeq
\noi Applying Lemma \ref{lemma:hhhh} with $h(\bbb{x})=f(\bbb{x})$, we conclude that $|\bbb{Ax}|_K=\bbb{0}$ will be achieved when $\rho \geq {L}/{\sigma(A)}$.

Since $|\bbb{Ax}|_K=\bbb{0}$, $\la (\bbb{Ax})_I, \bbb{u}_I \ra = \| (\bbb{Ax})_I\|_1$ and $\la (\bbb{Ax})_J, \bbb{u}_J \ra = \| (\bbb{Ax})_J\|_1$, we achieve that the complementarity constraint $\|\bbb{Ax}\|_1 - \la \bbb{Ax} , \bbb{u}\ra =0$ is fully satisfied. This finishes the proof for claim in the lemma.

\end{proof}
\end{lemma}

We now show that when the penalty parameter $\rho$ is larger than a threshold, the biconvex objective function $\J_{\rho}(\bbb{x},\bbb{u})$ is equivalent to the original constrained MPEC problem.

\begin{theorem} \label{theorem:1}
\textbf{Exactness of the Penalty Function.} The penalty problem $\min_{\bbb{x},\bbb{v}}~\J_{\rho}(\bbb{x},\bbb{v})$ admits the same local and global optimal solutions as the original MPEC problem when $\rho > L/\sigma(\bbb{A})$.
\begin{proof}

 First of all, based on the non-separable MPEC reformulation, we have the following Lagrangian function $\H: \mathbb{R}^n \times \mathbb{R}^m \times \mathbb{R} \rightarrow \mathbb{R}$:
\beq \label{eq:lag:nonsep}
\H(\bbb{x},\bbb{u},\chi) = f(\bbb{x}) + I(\bbb{u}) + \chi (\|\bbb{Ax}\|_1 - \la \bbb{Ax}, \bbb{u}\ra )\nn
\eeq
\noi Based on the Lagrangian function\footnote{Note that the Lagrangian function $\H$ has the same form as the penalty function $\J$ and the multiplier $\chi$ also plays a role of penalty parameter. It always holds at the optimal solution that $\chi^* > 0$ since $\la \bbb{Ax}, \bbb{u} \ra \leq \|\bbb{Ax} \|_1 \cdot \|\bbb{u}\|_{\infty}\leq \|\bbb{Ax} \|_1$ and the equality holds at the optimal solution.}, we have the following KKT conditions for any KKT solution ($\bbb{x}^*,\bbb{u}^*,\chi^*$):
\beq \label{eq:epm:kkt:0}
\zero & \in & \partial f(\bbb{x}^*)  + \chi^* \bbb{A}^T \partial | \bbb{Ax}^*| - \chi^* \bbb{A}^T\bbb{u}^* \nn\\
\zero & \in &    \partial I(\bbb{u}^{*}) - \chi^* \bbb{Ax}^*   \\
0 & = & \|\bbb{Ax}^*\|_1 - \la \bbb{Ax}^* , \bbb{u}^*\ra \nn
\eeq
\noi The KKT solution is defined as the first-order minimizer (respectively maximizer) for the primal (respectively dual) variables of the Lagrange function. It can be simply derived from setting the (sub-)gradient of the Lagrange function to zero with respect to each block of variables.

Secondly, we focus on the penalty function $\J_{\rho}(\bbb{x},\bbb{u})$. When $\rho > L/\sigma(\bbb{A})$, by Lemma \ref{lemma:v:01} we have:
\beq \label{eq:epm:kkt:2}
0=\|\bbb{Ax}\|_1 - \la \bbb{Ax} , \bbb{u}\ra
\eeq
\noi By the local optimality of the penalty function $\J_{\rho}(\bbb{x},\bbb{u})$, we have the following equalities:
\beq \label{eq:epm:kkt:1}
\begin{split}
\zero & \in & \partial f(\bbb{x})  + \rho \bbb{A}^T \partial | \bbb{Ax}| - \rho \bbb{A}^T\bbb{u}  \\
\zero & \in &    \partial I(\bbb{u}) - \rho \bbb{Ax}  ~~~~~~~~~~~~~~~~~~~~
\end{split}
\eeq
\noi Since Eq (\ref{eq:epm:kkt:1}) and Eq (\ref{eq:epm:kkt:2}) coincide with Eq (\ref{eq:epm:kkt:0}), we conclude that the solution of $\J_{\rho}(\bbb{x},\bbb{u})$ admits the same local and global optimal solutions as the original non-separable MPEC reformulation, when $\rho>L/\sigma(\bbb{A})$. 

\end{proof}
\end{theorem}

We now establish the convergence rate of the exact penalty method and determines the number of iterations beyond which a certain accuracy is guaranteed.

\begin{theorem}

\textbf{Convergence rate of Algorithm MPEC-EPM.} Assume that $\|\bbb{x}^t\|\leq \delta$ for all $t$. Algorithm 1 will converge to the first-order KKT point in at most $\lceil {\left(\ln(L \delta)-\ln(\epsilon \rho^0)\right)}/{ \ln2}\rceil$ outer iterations with the accuracy at least $ \|\bbb{Ax}\|_1 - \la \bbb{Ax},\bbb{v} \ra \leq \epsilon$. 
\begin{proof}

Assume Algorithm 1 takes $s$ outer iterations to converge. We have the following inequalities:
\beq
&& \textstyle \|\bbb{Ax}^{s+1}\|_1 - \la \bbb{Ax}^{s+1},\bbb{v}^{s+1} \ra \nn\\
&=& \textstyle \|(\bbb{Ax}^{s+1})_K\|_1 \nn\\
&\leq& \textstyle  \frac{ \delta L}{\rho^s} \nn
\eeq
\noi where the first step uses the notations and the results (See Eq (\ref{eq:KK})) in Lemma \ref{lemma:v:01}, the second step uses the third part of Lemma \ref{lemma:hhhh} with $h(\bbb{x})=f(\bbb{x})$. The above inequality implies that when $\rho^s\geq \frac{   { \delta  L}    }{\epsilon}$, Algorithm MPEC-EPM achieves accuracy at least $ \|\bbb{Ax}\|_1 - \la \bbb{Ax},\bbb{v} \ra \leq \epsilon$. Noticing that $\rho^s= 2^s \rho^0$, we have that $\epsilon$ accuracy will be achieved when
\beq
&& \textstyle 2^s \rho^0 \geq \frac{   { \delta L}    }{\epsilon }\nn\\
&\Rightarrow & \textstyle 2^s \geq \frac{   {\delta L}    }{\epsilon \rho^0 }\nn\\
&\Rightarrow & \textstyle s \geq{\left(\ln( \delta L ) -  \ln(\epsilon \rho^0)  \right)}/{ \ln2}       \nn
\eeq
\noi Thus, we finish the proof of this theorem.

\end{proof}
\end{theorem}

\section{ Alternating Direction Method} \label{sect:adm}
This section presents a proximal alternating direction method (PADM) for solving Eq (\ref{eq:l0}). This is mainly motivated by the recent popularity of ADM in the non-convex optimization literature. One direct solution is to apply ADM on the MPEC problem in Eq (\ref{eq:l0:mpec:nonsep}). However, this strategy may not be appealing, since it introduces a non-separable structure with no closed form solution for its sub-problems. Instead, we consider a separable MPEC reformulation used in \cite{YuanG15,BiLP14,yuan2016proximal}.

\subsection{Separable MPEC Reformulations}

 \begin{lemma}\label{lemma:mpec:sep}
(\cite{YuanG15}) For any given $\bbb{x}\in \mathbb{R}^{n}$, it holds that
  \beq\label{eq:mpec}
  \|\bbb{x}\|_0=\min_{\bbb{0}\le \bbb{v}\le \bbb{1}}~\la \bbb{1},\bbb{1}-\bbb{v}\ra,~s.t.~\bbb{v}\odot |\bbb{x}|=\bbb{0},
 \eeq
 and $\bbb{v}^*=\mathbf{1}-{\rm sign}(|\bbb{x}|)$ is the unique  solution of Eq(\ref{eq:mpec}).
 \end{lemma}
 \noi Using Lemma \ref{lemma:mpec:sep}, we can rewrite Eq (\ref{eq:l0}) in an equivalent form as follows.
\beq \label{eq:l0:mpec:sep}
\min_{\bbb{x},\bbb{v}}~f(\bbb{x})+I(\bbb{v}),~s.t.~|\bbb{Ax}|\odot \bbb{v}=\bbb{0}
\eeq
\noi where $I(\bbb{v})$ is an indicator function
on $\Omega$ with
\beq
I(\bbb{v}){\scriptsize=\begin{cases}
0,~\bbb{v} \in \Omega\\
\infty,~\bbb{v} \notin \Omega,\\\end{cases}}\Omega \triangleq \{\bbb{v}| \la \bbb{1}, \bbb{1}-\bbb{v} \ra \leq k, ~\bbb{0}\leq \bbb{v}\leq \bbb{1}\},~ \nn
\eeq

\begin{algorithm} [!h]
\caption{\label{alg:adm} { MPEC-ADM: An Alternating Direction Method for Solving MPEC Problem (\ref{eq:l0:mpec:sep})}}
\raggedright (S.0) Initialize $\bbb{x}^0=\bbb{0}\in\mathbb{R}^{n}$, $\bbb{v}^0 = \bbb{1} \in \mathbb{R}^{m}$, $\bbb{\pi}^{0}= \eta \cdot \bbb{1} \in \mathbb{R}^m$. Set $t=0$ and $\mu>0$.

\raggedright (S.1) Solve the following $\bbb{x}$-subproblem:
     \beq \label{eq:subprob:adm:x}
     \bbb{x}^{t+1}= \mathop{\arg\min}_{\bbb{x}}~\mathcal{L}(\bbb{x},\bbb{v}^{t},\bbb{\pi}^t) +  \frac{\mu}{2} \|\bbb{x}-\bbb{x}^t\|^2
     \eeq

\raggedright (S.2) Solve the following $\bbb{v}$-subproblem:
     \beq \label{eq:subprob:adm:v}
      \bbb{v}^{t+1} =\mathop{\arg\min}_{ \bbb{v}}~\mathcal{L} (\bbb{x}^{t+1},\bbb{v},\bbb{\pi}^t)+ \frac{\mu}{2} \|\bbb{v}-\bbb{v}^t\|^2
     \eeq

\raggedright (S.3) Update the Lagrange multiplier:
                \beq \label{eq:update:adm:pi}
                \bbb{\pi}^{t+1}=\bbb{\pi}^{t} + \alpha  (  |\bbb{Ax}^{t+1}| \odot \bbb{v}^{t+1}   )
                \eeq
\raggedright (S.4) Set $t:=t+1$ and then go to Step (S.1).
 \end{algorithm}

\subsection{Proposed Optimization Framework}



To solve Eq (\ref{eq:l0:mpec:sep}) using  PADM, we form the augmented Lagrangian function $\mathcal{L}(.)$
in Eq (\ref{eq:main:mpec:augmented:sep})
\beq \label{eq:main:mpec:augmented:sep}
\begin{split}
\mathcal{L}(\bbb{x},\bbb{v},\bbb{\pi}) \triangleq f(\bbb{x}) + I(\bbb{v}) +\la  |\bbb{Ax}|\odot \bbb{v} ,\bbb{\pi} \ra + \frac{\alpha}{2} \| |\bbb{Ax}|\odot \bbb{v} \|^2
\end{split}
\eeq
\noi where $\bbb{\pi}$ is the Lagrange multiplier associated with the constraint $|\bbb{Ax}|\odot \bbb{v}=\bbb{0}$, and $\alpha>0$ is the penalty parameter. We detail the PADM iteration steps for Eq (\ref{eq:l0:mpec:sep}) in Algorithm \ref{alg:adm}, which has the following properties. 

\vspace{1pt}\noi \bbb{(a) Initialization}. We set $\bbb{v}^0= \bbb{1}$ and $\bbb{\pi}^0=\eta \bbb{1}$, where $\eta$ is a small parameter. This finds a reasonable local minimum in the first iteration, as it reduces to an $\ell_1$ norm minimization problem for the $\bbb{x}$-subproblem. 

\vspace{1pt}\noi \bbb{(b) Monotone and boundedness property}. For any feasible solution $\bbb{v}$ in Eq (\ref{eq:subprob:adm:v}), it holds that $|\bbb{Ax}| \odot \bbb{v}\geq \bbb{0}$. Using the fact that $\alpha^t>0$ and due to the $\bbb{\pi}^t$ update rule, $\bbb{\pi}^t$ is monotone increasing. If we initialize $\bbb{\pi}^0>0$ in the first iteration, $\bbb{\pi}$ is always positive. Another key feature of this method is the boundedness of the multiplier $\bbb{\pi}$ (see Theorem \ref{theorem:3}). It reduces to alternating minimization algorithm for biconvex optimization problem \cite{Attouch2010,Bolte2014}.


\vspace{1pt}\noi \bbb{(c) $\bbb{v}$-Subproblem}. Variable $\bbb{v}$ in Eq (\ref{eq:subprob:adm:v}) is updated by solving the following problem:
\beq \label{eq:V:subproblem:1}
\begin{split}
 \bbb{v}^{t+1} =\mathop{\arg\min}_{\bbb{v}}   ~ \frac{1}{2} \bbb{v}^T\bbb{D}\bbb{v} + \la \bbb{v},\bbb{b}\ra,~s.t.~\bbb{v}\in \Omega
\end{split}
\eeq
\noi where $\bbb{b}\triangleq \bbb{\pi}^{t} \odot |\bbb{Ax}^{t+1}| - \mu \bbb{v}^t$ and $\bbb{D}$ is a diagonal matrix with $\bbb{d} \triangleq \alpha |\bbb{Ax}^{t+1}|\odot |\bbb{Ax}^{t+1}| + \mu$ in the main diagonal entries. Introducing the proximal term in the $\bbb{v}$-subproblem leads to a strongly convex problem. Thus, it can be solved \emph{exactly} using \cite{HelgasonKL80}.

\noi \bbb{(d) $\bbb{x}$-Subproblem}. Variable $\bbb{x}$ in Eq (\ref{eq:subprob:adm:x}) is updated by solving a re-weighted $\ell_1$ norm optimization problem. Similar to Eq (\ref{eq:subprob:epm:u}), it can be solved using classical/linearized ADM.


\subsection{Theoretical Analysis}


In the following theorem, we show that when the monotone increasing multiplier $\bbb{\pi}$ is larger than a threshold, the biconvex objective function in Eq (\ref{eq:main:mpec:augmented:sep}) is equivalent to the original constrained MPEC problem in Eq (\ref{eq:l0:mpec:sep}). Note that our proof is also built upon Lemma \ref{lemma:hhhh}.

\begin{theorem}  \label{theorem:3}
Any local optimal solution of the minimization problem: $\min_{\bbb{x},\bbb{v}}~\La(\bbb{x},\bbb{v},\bar{\bbb{\pi}})$ will be achieved with $|\bbb{Ax}| \odot \bbb{v} =\bbb{0}$, when $\bar{\bbb{\pi}}>L/\sigma(\bbb{A})$.




\begin{proof}
We assume that $\bbb{x}$ and $\bbb{v}$ are arbitrary local optimal solutions of the augmented Lagrangian function for a given $\bar{\bbb{\pi}}$.

\noi (i) Firstly, we now focus on the \bbb{v}-subproblem, we have:
\beq
\bbb{v} = \arg \min_{\bbb{v}}~I(\bbb{v})+\la\bar{\bbb{\pi}},|\bbb{Ax}|\odot \bbb{v} \ra + \frac{\alpha}{2} \||\bbb{Ax}|\odot \bbb{v}\|_2^2\nn
\eeq
\noi Then there exists a constant $\theta\geq0$ (that depends on the local optimal solutions $\bbb{x}$ and $\bbb{v}$) such that it solves the following minimax saddle point problem:
\beq
\max_{\theta\geq0}~\min_{\bbb{0}\leq \bbb{v}\leq \bbb{1}}~\la\bar{\bbb{\pi}},|\bbb{Ax}|\odot \bbb{v} \ra + \frac{\alpha}{2} \||\bbb{Ax}|\odot \bbb{v}\|_2^2 + \theta (n-k-\la \bbb{v},\bbb{1}\ra)\nn
\eeq
\noi Clearly, $\bbb{v}$ can be computed as:
\beq \label{eq:dual:v:bound}
\bbb{v} = \max(0,\min(1, \frac{\theta-|\bbb{Ax}|\odot \bar{\bbb{\pi}}}{\alpha |\bbb{Ax}|\odot |\bbb{Ax}|}))
\eeq
We now define $I\triangleq\{i|\theta-|\bbb{Ax}|_i \cdot \bar{\bbb{\pi}}_i \leq 0\}$,~$J\triangleq\{i|\theta-|\bbb{Ax}|_i \cdot \bar{\bbb{\pi}}_i > 0\}$. Clearly, we have
\beq
\bbb{v}_I=0,\bbb{v}_J\neq0\nn
\eeq
\noi (ii) Secondly, we focus on the $\bbb{x}$-subproblem:
\beq
\min_{\bbb{x}}~f(\bbb{x}) + \la  |\bbb{Ax}|\odot \bbb{v} ,\bbb{\pi} \ra + \frac{\alpha}{2} \| |\bbb{Ax}|\odot \bbb{v} \|^2\nn
\eeq

\noi Noticing $\bbb{v}_J \neq0$, we apply Lemma \ref{lemma:hhhh} with $\bbb{w}=\bar{\bbb{\pi}}\odot\bbb{v}$ and $h(\bbb{x})= f(\bbb{x}) + \frac{\alpha}{2} \| |\bbb{Ax}|\odot \bbb{v} \|^2$, then $|\bbb{Ax}|_{J}=0$ will be achieved whenever $\bar{\bbb{\pi}}_J\odot \bbb{v}_J>L_g/\sigma(\bbb{A})$, where $L_h$ is the Lipschitz constant of $h(\cdot)$. Since $\bbb{v}_I=0,~(\bbb{Ax})_J=0$ and $I \cup J=\{1,2,...,m\}$, we have $\frac{\alpha}{2} \||\bbb{Ax}|\odot \bbb{v} \|^2=0$, thus, $L_h= L$. Moreover, incorporating $(\bbb{Ax})_J$=0 into Eq(\ref{eq:dual:v:bound}), we have $\bbb{v}_J=\bbb{1}$.

\noi (iii) Finally, since $\bar{\bbb{\pi}}>0$ and $\theta\geq 0$, by the definition of $I$, we have $|\bbb{Ax}|_I \geq 0$. In summery, we have the following results:
\beq
\bbb{v}_I=0,~\bbb{v}_J=1,~(\bbb{Ax})_I \geq 0,~(\bbb{Ax})_J=0
\eeq
\noi We conclude that the complementarity constraint will hold automatically whenever $\bar{\bbb{\pi}}_J > ( L/\sigma(\bbb{A})) \odot \tfrac{1}{\bbb{v}_J} =  L/\sigma(\bbb{A})$, where $J$ is the index that $|\bbb{Ax}|_J=0$. However, we do not have any prior information of the index set $J$ (i.e. the value of $\bbb{Ax}^*$) and $J$ can be any subset of $\{1,2,...,m\}$, the condition $\bar{\bbb{\pi}}_J >   L/\sigma(\bbb{A})$ needs to be further restricted to $\bar{\bbb{\pi}} >  L/\sigma(\bbb{A})$. Thus, we complete the proof of this Lemma.




%
%
%
%
%

\end{proof}
\end{theorem}

\begin{theorem} \label{theorem:bound}
\textbf{Boundedness of multiplier and exactness of the Augmented Lagrangian Function.} The augmented Lagrangian problem $\min_{\bbb{x},\bbb{v}}~\La(\bbb{x},\bbb{v},\bar{\bbb{\pi}})$ admits the same local and global optimal solutions as the original MPEC problem when $\bar{\bbb{\pi}} \geq L/\sigma(\bbb{A})$.
\begin{proof}

First of all, based on the Lagrangian function $\La(\cdot)$, we have the following KKT conditions for any KKT solution ($\bbb{x}^*,\bbb{v}^*,\bbb{\pi}^*$):
\beq \label{eq:epm:kkt:00}
\zero & \in & \partial f(\bbb{x}^*)  + \bbb{A}^T (\partial | \bbb{Ax}^*| \odot \bbb{v}^* \odot \bbb{\pi}^*) \nn\\
\zero & \in &    \partial I(\bbb{v}^{*}) + \bbb{Ax}^* \odot \bbb{\pi}^*    \\
\bbb{0} & = & |\bbb{Ax}^*| \odot \bbb{v}^* \nn
\eeq
\noi The KKT solution is defined as the first-order minimizer (respectively maximizer) for the primal (respectively dual) variables of the Lagrange function. It can be simply derived from setting the (sub-)gradient of the Lagrange function to zero with respect to each block of variables.

Secondly, we focus on the penalty function $\La(\bbb{x},\bbb{u},\bar{\bbb{\pi}})$. When $\bar{\bbb{\pi}} > L/\sigma(\bbb{A})$, by Lemma \ref{lemma:equality:hold:2} we have:
\beq \label{eq:epm:kkt:3}
\bbb{0}= |\bbb{Ax}| \odot \bbb{v}
\eeq
\noi By the local optimality of the penalty function $\La(\bbb{x},\bbb{u},\bar{\bbb{\pi}})$, we have the following equalities:
\beq \label{eq:epm:kkt:4}
\zero & \in & \partial f(\bbb{x})  + \bbb{A}^T (\partial | \bbb{Ax}| \odot \bbb{v} \odot \bbb{\pi}) \nn\\
\zero & \in &    \partial I(\bbb{v}) + \bbb{Ax} \odot \bbb{\pi}
\eeq

\noi Since Eq (\ref{eq:epm:kkt:3}) and Eq (\ref{eq:epm:kkt:4}) coincide with Eq (\ref{eq:epm:kkt:00}), we conclude that the solution of $\La(\bbb{x},\bbb{u},\bar{\bbb{\pi}})$ admits the same local and global optimal solutions as the original non-separable MPEC reformulation, when $\bar{\bbb{\pi}}>L/\sigma(\bbb{A})$.

\end{proof}

\end{theorem}

In the following, we present the proof of Theorem 4. For the ease of discussions, we define:
\beq
\bbb{s} \triangleq \{\bbb{x},\bbb{v},\bbb{\pi}\},~\bbb{w} \triangleq \{\bbb{x},\bbb{v}\}.\nn
\eeq


\noi First of all, we prove the subgradient lower bound for the iterates gap by the following lemma.

\begin{lemma}\label{lemma:bound:grad}

Assume that $\bbb{x}^t$ are bounded for all $t$, then there exists a constant $\varpi>0$ such that the following inequality holds:
\beq \label{eq:bound:grad}
\| \partial \La (\bbb{s}^{t+1}) \| \leq \varpi  \| \bbb{s}^{t+1} - \bbb{s}^{t}\|
\eeq
\begin{proof}
For notation simplicity, we denote $\bbb{p} = \partial |\bbb{Ax}^{t+1}|$ and $\bbb{z} =|\bbb{Ax}^{t+1}|$. By the optimal condition of the $\bbb{x}$-subproblem and $\bbb{v}$-subproblem, we have:
{\beq \label{eq:opt:X}
\zero &\in&  \bbb{A}^T \left( \bbb{p} \odot \bbb{v}^{t} \odot ( \alpha \bbb{z} \odot \bbb{v}^{t} + \bbb{\pi}^{t})\right) +  f'(\bbb{x}^{t+1})+  \mu  (\bbb{x}^{t+1}-\bbb{x}^t) \nn\\
\zero &\in&   \left( \alpha \bbb{z}  \odot \bbb{v}^{t+1} + \bbb{\pi}^{t}\right) \odot \bbb{z} +  \partial I(\bbb{v}^{t+1}) + \mu(\bbb{v}^{t+1}-\bbb{v}^t)
\eeq}
\noi By the definition of $\La$ we have that
{\small
\beq \label{eq:bound:Lx}
&&\partial \La_{\bbb{x}} (\bbb{s}^{t+1})\nn \\
&=&    \bbb{A}^T (\bbb{p} \odot \bbb{v}^{t+1} \odot \alpha ( \bbb{z}  \odot \bbb{v}^{t+1} + \bbb{\pi}^{t+1}/\alpha)) +  f'(\bbb{x}^{t+1})     \nn \\
&=&  \bbb{A}^T \bbb{p} \odot ( \bbb{v}^{t+1} \odot \alpha ( \bbb{z} \odot \bbb{v}^{t+1} + \bbb{\pi}^{t+1}/\alpha) \nn\\
&&-  \bbb{v}^{t} \odot \alpha (\bbb{z}  \odot \bbb{v}^{t} + \bbb{\pi}^{t}/\alpha)  )     + \mu(\bbb{x}^{t}-\bbb{x}^{t+1}) \nn \\
&=&  \bbb{A}^T \bbb{p} \odot ( \bbb{v}^{t+1} \odot \alpha ( \bbb{z} \odot \bbb{v}^{t+1} + \bbb{\pi}^{t+1}/\alpha  \nn\\
&& -  \bbb{v}^{t} \odot \alpha ( \bbb{z}  \odot \bbb{v}^{t} + (\bbb{\pi}^{t+1} - \bbb{\pi}^{t+1} + \bbb{\pi}^{t})/\alpha)  )   + \mu(\bbb{x}^{t}-\bbb{x}^{t+1}) \nn \\
&=& \bbb{A}^T \bbb{p} \odot \bbb{v}^t \odot (\bbb{\pi}^{t} - \bbb{\pi}^{t+1}) + \mu(\bbb{x}^{t}-\bbb{x}^{t+1}) \nn\\
 &&+ \bbb{A}^T \bbb{p} \odot ( \bbb{v}^{t+1} \odot  ( \alpha \bbb{z}  \odot \bbb{v}^{t+1} + \bbb{\pi}^{t+1})  -  \bbb{v}^{t} \odot  ( \alpha \bbb{z}  \odot \bbb{v}^{t} + \bbb{\pi}^{t+1} )  )    \nn \\
&=& \bbb{A}^T \bbb{p} \odot \bbb{v}^t \odot (\bbb{\pi}^{t} - \bbb{\pi}^{t+1}) + \mu(\bbb{x}^{t}-\bbb{x}^{t+1}) \nn\\
&&   +   \bbb{A}^T \bbb{p} \odot \bbb{\pi}^{t+1} \odot ( \bbb{v}^{t+1} - \bbb{v}^{t} )  +  \bbb{A}^T \bbb{p} \odot   \alpha \bbb{z}  \odot  ( \bbb{v}^{t+1} \odot  \bbb{v}^{t+1}  - \bbb{v}^{t} \odot     \bbb{v}^{t}   )
\eeq
}\noi The first step uses the definition of $\La_{\bbb{x}} (\bbb{s}^{t+1})$, the second step uses Eq (\ref{eq:opt:X}), the third step uses $\bbb{v}^t + \bbb{v}^{t+1} - \bbb{v}^{t+1} = \bbb{v}^t$, the fourth step uses the multiplier update rule for $\bbb{\pi}$.

\noi We assume that $\bbb{x}^{t+1}$ is bounded by $\delta$ for all $t$, i.e. $\|\bbb{x}^{t+1}\|\leq \delta$. By Theorem \ref{theorem:bound}, $\bbb{\pi}^{t+1}$ is also bounded. We assume it is bounded by a constant $\kappa$, i.e. $\|\bbb{\pi}^{t+1}\|\leq \kappa$. Using the fact that $\|\partial \|\bbb{Ax}\|_1\|_{\infty} \leq 1$,~$\|\bbb{Ax}\|_{2} \leq \|\bbb{A}\|\|\bbb{x}\| $,~$\|\bbb{x}\odot \bbb{y}\|_{2} \leq \|\bbb{x}\|_{\infty} \|\bbb{y}\| $ and Eq (\ref{eq:bound:Lx}), we have the following inequalities:
{\small
\beq \label{eq:subg:bound:X}
&&\| \partial \La_{\bbb{x}} (\bbb{s}^{t+1})\|\nn\\
&\leq &  \|\bbb{A}\| \cdot \| \bbb{p}\|_{\infty} \cdot \|\bbb{v}\|_{\infty} \cdot \|\bbb{\pi}^{t} - \bbb{\pi}^{t+1}\| + \mu \|\bbb{x}^{t}-\bbb{x}^{t+1}\| + \nn \\
&& \|\bbb{A}\| \cdot \| \bbb{p}\|_{\infty} \cdot \|\bbb{\pi}^{t+1} \|_{\infty} \cdot \| \bbb{v}^{t+1} - \bbb{v}^{t} \| + \nn\\
&& \|\bbb{A}\| \cdot \| \bbb{p}\|_{\infty} \cdot \alpha \|\bbb{z}\| \cdot \|\bbb{v}^{t+1}  + \bbb{v}^{t} \|_{\infty} \cdot  \| \bbb{v}^{t+1}  - \bbb{v}^{t} \| \nn\\
&\leq& \|\bbb{A}\| \|\bbb{\pi}^{t} - \bbb{\pi}^{t+1}\| + \mu \|\bbb{x}^{t+1} -\bbb{x}^t\| +  (\kappa+2 \delta \alpha \|\bbb{A}\| )\cdot \|\bbb{A}\| \cdot \|\bbb{v}^{t+1}-\bbb{v}^t\|
\eeq
}
\noi Similarly, we have
\beq
 \textstyle \partial \La_{\bbb{v}}(\bbb{s}^{t+1})&=& \textstyle \partial I (\bbb{v}^{t+1}) + \alpha \bbb{z}\odot ( \bbb{\pi}^{t+1} /\alpha +  \bbb{z} \odot \bbb{v}^{t+1}) \nn \\
&=& \textstyle \bbb{z} \odot ( \bbb{\pi}^{t+1} - \bbb{\pi}^{t} )  - \mu (\bbb{v}^{t+1}-\bbb{v}^t) \nn
\eeq
\beq
\textstyle \partial \La_{\pi} (\bbb{s}^{t+1}) =|\bbb{Ax}^{t+1} \odot \bbb{v}^{t+1} |= \frac{1}{\alpha}\nn (\bbb{\pi}^{t+1} - \bbb{\pi}^{t})
\eeq
\noi Then we derive the following inequalities:
\beq
\textstyle \|\partial \La_{\bbb{v}}(\bbb{s}^{t+1}) \| \leq \delta \|\bbb{A}\| \|\bbb{\pi}^t -  \bbb{\pi}^{t+1}\| + \mu  \|\bbb{v}^{t+1}-\bbb{v}^t\| \label{eq:subg:bound:V}
\eeq
\vspace{-20pt}
\beq
\textstyle \|\partial \La_{\pi}(\bbb{s}^{t+1}) \| \leq \frac{1}{\alpha} \|\bbb{\pi}^{t+1} - \bbb{\pi}^{t}\| ~~~~~~~~~~~~~~~~~~~~~~~~~ \label{eq:subg:bound:pi}
\eeq
\noi Combining Eqs (\ref{eq:subg:bound:X}-\ref{eq:subg:bound:pi}), we conclude that there exists $\varpi>0$ such that the following inequality holds
\beq
\textstyle \| \partial \La(s^{t+1}) \| \leq \varpi \|\bbb{s}^{t+1}\| \nn
\eeq
\noi We thus complete the proof of this lemma.
\end{proof}
\end{lemma}

The following lemma is useful in our convergence analysis.

\begin{proposition} \label{lemma:cluster:piint}
Assume that $\bbb{x}^t$ are bounded for all $t$, then we have the following inequality:
\beq
\textstyle \sum_{t=0}^{\infty} \| \bbb{s}^t - \bbb{s}^{t+1} \|^2 < +\infty\nn
\eeq
\noi In particular the sequence $\|\bbb{s}^t-\bbb{s}^{t+1}\|$ is asymptotic regular, namely $\|\bbb{s}^t-\bbb{s}^{t+1}\|\rightarrow0$ as $t\rightarrow\infty$. Moreover any cluster point of $\bbb{s}^t$ is a stationary point of $\La$.

\begin{proof}
Due to the initialization and the update rule of $\bbb{\pi}$, we conclude that $\bbb{\pi}^t$ is nonnegative and monotone non-decreasing. Moreover, using the result of Theorem \ref{theorem:bound}, as $t\rightarrow \infty$, we have: $|\bbb{Ax}^{t+1}|\odot \bbb{v}^{t+1} = \bbb{0}$. Therefore, we conclude that as $t \rightarrow +\infty$ it must hold that
\beq \label{eq:bound:pi}
\textstyle |\bbb{Ax}^{t+1}|\odot \bbb{v}^{t+1} = \bbb{0},\nn\\
\textstyle \sum_{i=1}^t \|\bbb{\pi}^{i+1} - \bbb{\pi}^i\| < +\infty,\nn\\
\textstyle \sum_{i=1}^t \|\bbb{\pi}^{i+1} - \bbb{\pi}^i\|^2 < +\infty\nn
\eeq
\noi On the other hand, we naturally derive the following inequalities:
{\small\beq
&&\textstyle \La(\bbb{x}^{t+1},\bbb{v}^{t+1};\bbb{\pi}^{t+1})\nn\\
&=& \textstyle \La(\bbb{x}^{t+1},\bbb{v}^{t+1};\bbb{\pi}^{t}) + \la \bbb{\pi}^{t+1}-\bbb{\pi}^{t} , |\bbb{Ax}^{t+1}|\odot \bbb{v}^{t+1} \ra   \nn \\
&=& \textstyle \La(\bbb{x}^{t+1},\bbb{v}^{t+1};\bbb{\pi}^{t}) + \frac{1}{\alpha} \|\bbb{\pi}^{t+1} - \bbb{\pi}^t\|^2  \nn \\
&\leq& \textstyle \La(\bbb{x}^t,\bbb{v}^{t+1};\bbb{\pi}^t) - \frac{\mu}{2} \|\bbb{x}^{t+1} - \bbb{x}^t\|^2  + \frac{1}{\alpha} \|\bbb{\pi}^{t+1} - \bbb{\pi}^t\|^2 \label{eq:lowbound:L1}\nn \\
&\leq& \textstyle \La(\bbb{x}^{t},\bbb{v}^{t};\bbb{\pi}^{t})  - \frac{\mu}{2} \|\bbb{x}^{t+1} - \bbb{x}^t\|^2 - \frac{\mu}{2} \|\bbb{v}^{t+1} - \bbb{v}^t\|^2+ \textstyle \frac{1}{\alpha} \|\bbb{\pi}^{t+1} - \bbb{\pi}^t\|^2 \label{eq:lowbound:L}
\eeq}
\noi The first step uses the definition of $\La$; the second step uses update rule of the Lagrangian multiplier $\bbb{\pi}$; the third and fourth step use the $\mu$-strongly convexity of $\La$ with respect to $\bbb{x}$ and $\bbb{v}$, respectively.

We define $C=\frac{1}{\alpha} \sum_{i=1}^t \|\bbb{\pi}^{i+1} - \bbb{\pi}^i\|^2 + \La(\bbb{x}^{0},\bbb{v}^{0};\bbb{\pi}^0) - \La(\bbb{x}^{t+1},\bbb{v}^{t+1};\bbb{\pi}^{t+1})$. Clearly, by the boundedness of $\{\bbb{x}^t,\bbb{v}^t,\bbb{\pi}^t\}$, both $C$ and $\La(\bbb{x}^{t+1},\bbb{v}^{t+1};\bbb{\pi}^{t+1})$ are bounded. Summing Eq (\ref{eq:lowbound:L}) over $i=1, 2..., t$, we have:
\beq \label{eq:bound:w}
\textstyle\frac{\mu}{2} \sum_{i=1}^{t}~\|\bbb{w}^{i+1} - \bbb{w}^i\|^2  &\leq& C
\eeq
\noi Therefore, combining Eq(\ref{eq:bound:pi}) and Eq (\ref{eq:bound:w}), we have $\sum_{t=1}^{+\infty} \|\bbb{s}^{t+1} - \bbb{s}^{t}\|^2 < +\infty$; in particular $\| \bbb{s}^{t+1} - \bbb{s}^{t} \|\rightarrow 0$. By Eq (\ref{eq:bound:grad}), we have that:
\beq
 \|\partial \La(\bbb{s}^{t+1})\| \leq \varpi \|\bbb{s}^{t+1} - \bbb{s}^{t}\| \rightarrow 0\nn
\eeq
\noi which implies that any cluster point of $\bbb{s}^t$ is a stationary point of $\La$. We complete the proof of this lemma.

\end{proof}
\end{proposition}

\noi \textbf{Remarks}: Lemma \ref{lemma:cluster:piint} states that any cluster point is the KKT point. Strictly speaking, this result does not imply the convergence of the algorithm. This is because the boundedness of $\sum_{t=0}^{\infty} \| \bbb{s}^t - \bbb{s}^{t+1} \|^2$ does not imply that the sequence $\bbb{s}^t$ is convergent \footnote{One typical counter-example is $\bbb{s}^t=\sum_{i=1}^t \frac{1}{i}$. Clearly, $\sum_{t=0}^{\infty} \| \bbb{s}^t - \bbb{s}^{t+1} \|^2 = \sum_{t=1}^{\infty} (\frac{1}{t})^2$ is bounded by $\frac{\pi^2}{6}$; however, $\bbb{s}^t$ is divergent since $s^t=\ln(k)+C_{\text{e}}$, where $C_{\text{e}}$ is the well-known Euler's constant.}. In what follows, we aim to prove stronger result in Theorem \ref{theorem:adm:main}.

Our analysis is mainly based on a recent non-convex analysis tool called Kurdyka-{\L}ojasiewicz inequality \cite{Attouch2010,Bolte2014}. One key condition of our proof requires that the Lagrangian function $\La(\bbb{s})$ satisfies the so-call (KL) property in its effective domain. It is so-called the semi-algebraic function satisfy the Kurdyka-{\L}ojasiewicz property. We note that semi-algebraic functions include (i) real polynomial functions, (ii) finite sums and products of semi-algebraic functions , and (iii) indicator functions of semi-algebraic sets \cite{Bolte2014}. Using these definitions repeatedly, the graph of $\La(\bbb{s}): \{(\bbb{s},z)|z=\La(\bbb{s})\}$ can be proved to be a semi-algebraic set. Therefore, the Lagrangian function $\La(\bbb{s})$ is a semi-algebraic function. This is not surprising since semi-algebraic function is ubiquitous in applications \cite{Bolte2014}. We now present the following proposition established in \cite{Attouch2010}.


%
%

\begin{proposition}\label{prop:dist ineq}
For a given semi-algebraic function $\La(s)$, for all $s \in \text{dom} \La$, there exists $\theta\in[0,1),~\eta \in (0,+\infty]$ a neighborhood $\mathcal{S}$ of $s$ and a concave and continuous function $\varphi(\varsigma)= \varsigma^{1-\theta}$, $\varsigma\in [0,\eta)$ such that for all $\bar{\bbb{s}}\in \mathcal{S}$ and satisfies $\La(\bar{s}) \in (\La(s),\La(s)+\eta)$, the following inequality holds:
\begin{align}
dist(0,\partial \La(\bar{\bbb{s}})) \varphi'( \La(\bbb{s}) - \La(\bar{\bbb{s}}) ) \geq 1,~\forall \bbb{s}\nn \label{eq:kl:inequality}
\end{align}
\noi where $dist(0,\partial \La(\bar{\bbb{s}}))=\min \{||w^*||: w^* \in  \partial \La(\bar{\bbb{s}})  \}$.

\end{proposition}



%

Based on the Kurdyka-{\L}ojasiewicz inequality \cite{Attouch2010,Bolte2014}, the following theorem establishes the convergence of the alternating direction method.

\begin{theorem}\label{theorem:adm:main}
Assume that $\bbb{x}^t$ are bounded for all $t$. Then we have the following inequality:
\beq
\textstyle \sum_{t=0}^{+\infty} \|\bbb{s}^t-\bbb{s}^{t+1}\| < \infty
\eeq
Moreover, as $t \rightarrow +\infty$, Algorithm MPEC-ADM converges to a first order KKT point of the reformulated MPEC problem.
\begin{proof}

For simplicity, we define $r^t=\varphi( \La(\bbb{s}^{t}) - \La(\bbb{s}^{*}) )  - \varphi( \La(\bbb{s}^{t+1}) - \La(\bbb{s}^{*}))$. We naturally derive the following inequalities:
{\small
\beq \label{eq:y:k1k}
&&\textstyle \frac{\mu}{2} \|\bbb{w}^{t+1} - \bbb{w}^t\|^2  - \frac{1}{\alpha} \|\bbb{\pi}^{t+1} - \bbb{\pi}^t\|^2 \nn\\
&\leq& \textstyle\La(\bbb{s}^{t}) - \La(\bbb{s}^{t+1}) \nn\\
&=& \textstyle (\La(\bbb{s}^{t}) - \La(\bbb{s}^{*})) - (\La(\bbb{s}^{t+1})  - \La(\bbb{s}^{*})) \nn\\
&\leq& \textstyle  {r^t}/{  \varphi'( \La(\bbb{s}^{t}) - \La(\bbb{s}^{*}) )  } \nn \\
&\leq& \textstyle r^t  \cdot dist(0,\partial \La(\bbb{s}^{t})) \nn \\
&\leq& \textstyle r^t \cdot \varpi \left( \| \bbb{x}^{t} - \bbb{x}^{k-1}\| + \|\bbb{v}^{t} - \bbb{v}^{k-1}\| + \|\bbb{\pi}^{t} - \bbb{\pi}^{k-1}\| \right)\nn\\
&\leq& \textstyle  r^t \cdot \varpi \left( \sqrt{2} \| \bbb{w}^{t} - \bbb{w}^{k-1}\|  + \|\bbb{\pi}^{t} - \bbb{\pi}^{k-1}\| \right)\nn\\
&=& \textstyle  r^t \varpi \sqrt{\frac{16}{\mu}} \left(  \sqrt{\frac{\mu}{8}} \| \bbb{w}^{t} - \bbb{w}^{k-1}\|  +  \sqrt{\frac{\mu}{16}} \|\bbb{\pi}^{t} - \bbb{\pi}^{k-1}\| \right)\nn\\
&\leq& \textstyle   (  \sqrt{\frac{\mu}{8}} \| \bbb{w}^{t} - \bbb{w}^{k-1}\|  +  \sqrt{\frac{\mu}{16}} \|\bbb{\pi}^{t} - \bbb{\pi}^{k-1}\| )^2 + (r^t \varpi \sqrt{\frac{4}{\mu}})^2 \nn
\eeq
}
\noi The first step uses Eq (\ref{eq:lowbound:L}); the third step uses the concavity of $\varphi$ such that $\varphi(a)-\varphi(b) \geq \varphi'(a) (a-b)$ for all $a,b\in \mathbb{R}$; the fourth step uses the KL property such that $dist(0,\partial \La(\bbb{s}^{t})) \varphi'( \La(\bbb{s}^{t}) - \La(\bbb{s}^{*}) ) \geq 1$ as in Proposition \ref{prop:dist ineq}; the fifth step uses Eq (\ref{eq:bound:grad}); the sixth step uses the inequality that $\| x;y\| \leq  \|x\|+\|y\| \leq \sqrt{2}\| x;y\|$, where `;' in $[\cdot]$ denotes the row-wise partitioning indicator as in Matlab; the last step uses that fact that $2ab\leq a^2+b^2$ for all $a,b\in \mathbb{R}$.
\noi As $k\rightarrow \infty$, we have:
\beq
\textstyle  \frac{\mu}{2} \|\bbb{w}^{t+1} - \bbb{w}^t\|^2  \leq   \frac{4(r^t \kappa)^2}{\mu} + \left(  \sqrt{\frac{\mu}{8}} \| \bbb{w}^{t} - \bbb{w}^{k-1}\| \right)^2\nn
\eeq
\noi Taking the squared root of both side and using the inequality that $\sqrt{a+b}\leq \sqrt{a}+\sqrt{b}$ for all $a,b\in \mathbb{R}_+$, we have
\beq
\textstyle 2\sqrt{\frac{\mu}{8}} \|\bbb{w}^{t+1} - \bbb{w}^t\|  \leq   \frac{2r^t \kappa}{\sqrt{\mu}} +  \sqrt{\frac{\mu}{8}} \| \bbb{w}^{t} - \bbb{w}^{k-1}\|\nn
\eeq
\noi Then we have
{\small
\beq \label{eq:y:k1k}
\textstyle \sqrt{\frac{\mu}{8}} \|\bbb{w}^{t+1} - \bbb{w}^t\| \label{eq:y:k1k} \leq \frac{2r^t \kappa}{\sqrt{\mu}} +  \sqrt{\frac{\mu}{8}} (\| \bbb{w}^{t} - \bbb{w}^{k-1}\| - \|\bbb{w}^{t+1} - \bbb{w}^t\|)
\eeq}
\noi Summing Eq (\ref{eq:y:k1k}) over $i=1,2..., t$, we have:
\beq \label{eq:bound:0}
\textstyle \sum_{i=1}^t \sqrt{\frac{\mu}{8}}  \|\bbb{w}^{i+1} - \bbb{w}^i\| \leq \frac{\kappa}{\sqrt{\mu}}\sum_{i=1}^t  r^i  + \left( \| \bbb{w}^{1} - \bbb{w}^{0}\|  + \|\bbb{w}^{t+1} - \bbb{w}^{t}\| \right)~~~~~~~
\eeq
\noi The first term in the right-hand side of Eq(\ref{eq:bound:0}) is bounded since $\sum_{i=1}^t  r^i = \varphi( \La(\bbb{s}^{0}) - \La(\bbb{s}^{*}) )  - \varphi( \La(\bbb{s}^{t+1}) - \La(\bbb{s}^{*}))$ is bounded. Therefore, we conclude that as $k\rightarrow \infty$, we obtain:
\beq
\textstyle \sum_{i=1}^t   \|\bbb{w}^{i+1} - \bbb{w}^i\| < + \infty\nn
\eeq
By the boundedness of $\bbb{\pi}$ in Theorem \ref{theorem:bound}, we have $\textstyle \sum_{i=1}^t   \|\bbb{\pi}^{i+1} - \bbb{\pi}^i\| < + \infty$. Therefore, we obtain:
\beq
\textstyle \sum_{i=1}^t   \|\bbb{s}^{i+1} - \bbb{s}^i\| < + \infty\nn
\eeq
Finally, by Eq (\ref{eq:bound:grad}) in Lemma \ref{lemma:bound:grad} we have $\partial \La(\bbb{s}^{t+1})=0$. In other words, we have the following results:
\beq
\zero & \in &  f'(\bbb{x}^{t+1})  + \bbb{A}^T (\partial | \bbb{Ax}^{t+1}| \odot \bbb{v}^{t+1} \odot \bbb{\pi}^{t+1})   \nn\\
\zero & \in & \bbb{\pi}^{t+1} \odot |\bbb{Ax}^{t+1}|   + \partial I(\bbb{v}^{t+1})   \nn\\
\zero & = & |\bbb{Ax}^{t+1}|  \odot \bbb{v}^{t+1} \nn
\eeq
\noi which coincide with the KKT conditions in Eq (\ref{eq:epm:kkt:00}). Therefore, we conclude that the solution $\{ \bbb{x}^{t+1}, \bbb{v}^{t+1}, \bbb{\pi}^{t+1} \}$ converges to a first-order KKT point of the reformulated MPEC problem.

\end{proof}
\end{theorem}

\noi In above theorem, we assume that the solution $\bbb{x}$ is bounded. This assumption is inherited from the nonconvex and nonsmooth minimization algorithm (see Theorem 1 in \cite{Bolte2014}, where they make the same assumption). In fact, with minor modification to the proofs, the boundedness assumption can be removed by adding a compact set constraint $\|\bbb{x}\|_{\infty} < \delta$ to Eq (\ref{eq:l0}), for the case $\bbb{A}=\bbb{I}$.



\section{On MPEC Optimization}\label{sect:disc}

In this paper, MPEC reformulations are considered to solve the $\ell_0$ norm problem. Mathematical programs with equilibrium constraints (MPEC) are optimization problems where the constraints include variational inequalities or complementarities. It is related to the Stackelberg game and it is used in the study of economic equilibrium and engineering design. MPECs are difficult to deal with because their feasible region is not necessarily convex or even connected.

\textbf{General MPEC motivation.} The basic idea behind MPEC-based methods is to transform the thin and nonsmooth, nonconvex feasible region into a thick and smooth one by introducing regularization/penalization on the complementary/equilibrium constraints. In fact, several other nonlinear methods solve the MPEC problem with similar motivations. For example, an exact penalty method where the complementarity term is moved to the objective in the form of an $\ell_1$-penalty was proposed in \cite{hu2004convergence}. The work in \cite{facchinei1999smoothing} suggested a smoothing family by replacing the complementarity system with the perturbed Fischer-Burmeister function. Also, a quadratic regularization technique to handle the complementarity constraints was proposed in \cite{demiguel2005two}. 

\begin{figure} [!h]
\begin{center}
\includegraphics[width=0.55\textwidth,height=0.45\textwidth]{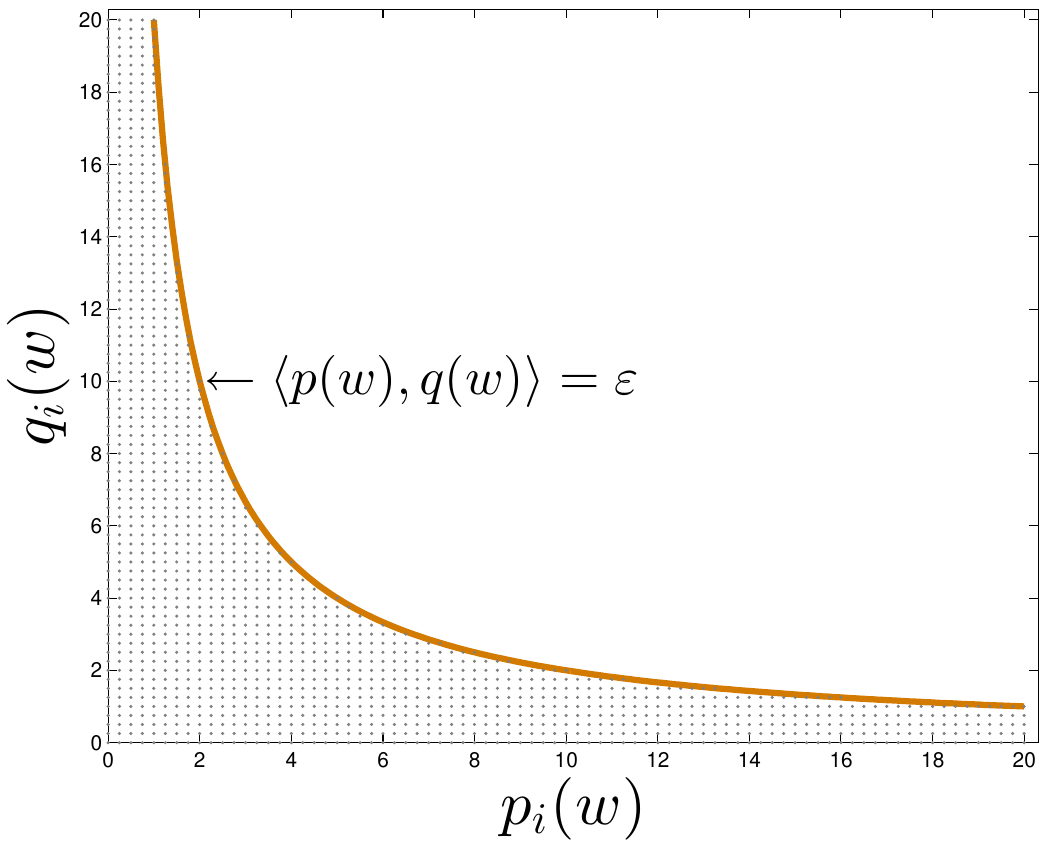}
\vspace{1pt}\caption{Effect of regularizing/smoothing the original optimization problem. }
\label{fig:mpec}
\end{center}
\end{figure}

\textbf{In the defense of our methods.} We propose two different penalization/regularization schemes to solve the equivalent MPECs. They have several merits. \textbf{(a)} Each reformulation is equivalent to the original $\ell_0$ norm problem. \textbf{(b)} They are continuous reformulations, so they facilitate KKT analysis and are amenable to the use of existing continuous optimization techniques to solve the convex sub-problems. We argue that, from a practical viewpoint, improved solutions to Eq (\ref{eq:l0}) can be obtained by reformulating the problem using MPEC and focusing on the complementarity constraints. \textbf{(c)} They find a good initialization because they both reduce to a convex relaxation method in the first iteration. \textbf{(d)} Both methods exhibit strong convergence guarantees, this is because they essentially reduce to alternating minimization for bi-convex optimization problems \cite{Attouch2010,Bolte2014}. \textbf{(e)} They have a monotone/greedy property owing to the complimentarity constraints brought on by MPEC. The complimentary system characterizes the optimality of the KKT solution. We let $w \triangleq \{\bbb{x},\bbb{u}\}$ (or $w \triangleq \{\bbb{x},\bbb{v}\}$). Our solution directly handles the complimentary system of Eq (\ref{eq:l0}) which takes the following form: $\la p(w),q(w) \ra = 0,~p(w) \geq  0,~ q(w) \geq 0$. Here $p(\cdot)$ and $q(\cdot)$ are non-negative mappings of $w$ which are vertical to each other. The complimentary constraints enable all the special properties of MPEC that distinguish it from general nonlinear optimization. We penalize the complimentary error of $\varepsilon \triangleq \la p(w),q(w)\ra$ (which is always non-negative) and ensure that the error $\varepsilon$ is decreasing in every iteration. See Figure \ref{fig:mpec} for a geometric interpretation for MPEC optimization.

\textbf{MPEC-EPM vs. MPEC-ADM.} We consider two algorithms based on two variational characterizations of the $\ell_0$ norm function (separable and non-separable MPEC\footnote{Note that the non-separable MPEC has one equilibrium constraint ($\|\bbb{x}\|_1= \la \bbb{x},\bbb{u}\ra$) and the separable MPEC has $n$ equilibrium constraints ($|\bbb{x}|\odot \bbb{v}=\bbb{0}$). The terms separable and non-separable are related to whether the constraints can be decomposed to independent components.}). Generally speaking, both methods have their own merits. \textbf{(a)} MPEC-EPM is more simple since it can directly use existing $\ell_1$ norm minimization solvers/codes while MPEC-ADM involves additional computation for the quadratic term in the standard Lagrangian function. \textbf{(b)} MPEC-EPM is less adaptive in its per-iteration optimization, since the penalty parameter $\rho$ is monolithically increased until a threshold is achieved. In comparison, MPEC-ADM is more adaptive, since a constant penalty also guarantees monotonically non-decreasing multipliers and convergence. \textbf{(c)} Regarding numerical robustness, while MPEC-EPM needs to solve the subproblem with certain accuracy, MPEC-ADM is more robust since the dual Lagrangian function provides a support function and the solution never degenerates even if the $\bbb{x}$-problem is solved only approximately. 

\section{Experimental Validation} \label{sect:exp}
\vspace{-8pt}
In this section, we demonstrate the effectiveness of our algorithms on five $\ell_0$ norm optimization tasks, namely feature selection, segmented regression, trend filtering, MRF optimization, and impulse noise removal. We compare our proposed MPEC-EPM (Algorithm \ref{alg:epm}) and MPEC-ADM (Algorithm \ref{alg:adm}) against the following methods.
\begin{itemize}
  \item Greedy Hard Thresholding (GREEDY) \cite{beck2013sparsity}. It considers decreasing the objective function and identifying the active variables simultaneously using the following update: $\bbb{x}  \Leftarrow \arg \min_{\|\bbb{y}\|_0\leq k} \|\bbb{x}-f'(\bbb{x})/L_f-\bbb{y}\|_2^2$, where $L_f$ is the gradient Lipschitz continuity constant. Note that this method is not applicable when the objective is non-smooth.

  \item Gradient Support Pursuit (GSP). This is a state-of-the-art greedy algorithm that approximates sparse minima of cost functions which have stable restricted Hessian. We use the Matlab implementation provided by the author \url{users.ece.gatech.edu/sbahmani7/GraSP.html}.

  \item Convex $\ell_1$ minimization method (CVX). Since this approximate method can not control the sparsity of the solution, we solve a $\ell_1$ regularized problem where the regulation parameter is swept over $\{2^{-10},2^{-8},...,2^{10}\}$. Finally, the solution that leads to smallest objective value after a hard thresholding projection (which reduces to setting the smallest $n-k$ values of the solution in magnitude to 0) is selected.

  \item Quadratic Penalty Method (QPM) \cite{LuZ13}. This is splitting method that introduces auxiliary variables to separate the calculation of the non-differentiable and differentiable term and then performs block coordinate descend on each subproblem. We use GREEDY to solve the subproblem with respect to $\bbb{x}$ with reasonable accuracy. We increase the penalty parameter by $\sqrt{10}$ in every $T=30$ iterations.

  \item Direct Alternating Direction Method (DI-ADM). This is an ADM directly applied to Eq (\ref{eq:l0}). We use a similar experimental setting as QPM.

  \item Mean Doubly Alternating Direction Method (MD-ADM) \cite{Dong2013}. This is an ADM that treats arithmetic means of the solution sequence as the actual output. We also use a similar parameter setting as QPM.

\end{itemize}

\begin{figure*}[!t]
\captionsetup[subfigure]{justification=centering}
    \centering
      \begin{subfigure}{\imgwidlogloss}\includegraphics[width=\textwidth,height=\imgheilogloss]{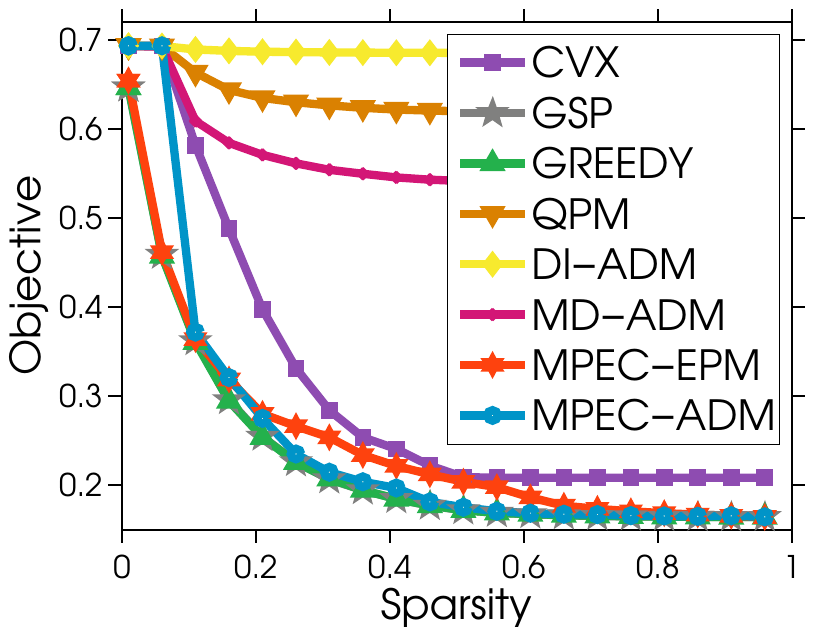}\caption{`a6a'}\end{subfigure}\ghs
      \begin{subfigure}{\imgwidlogloss}\includegraphics[width=\textwidth,height=\imgheilogloss]{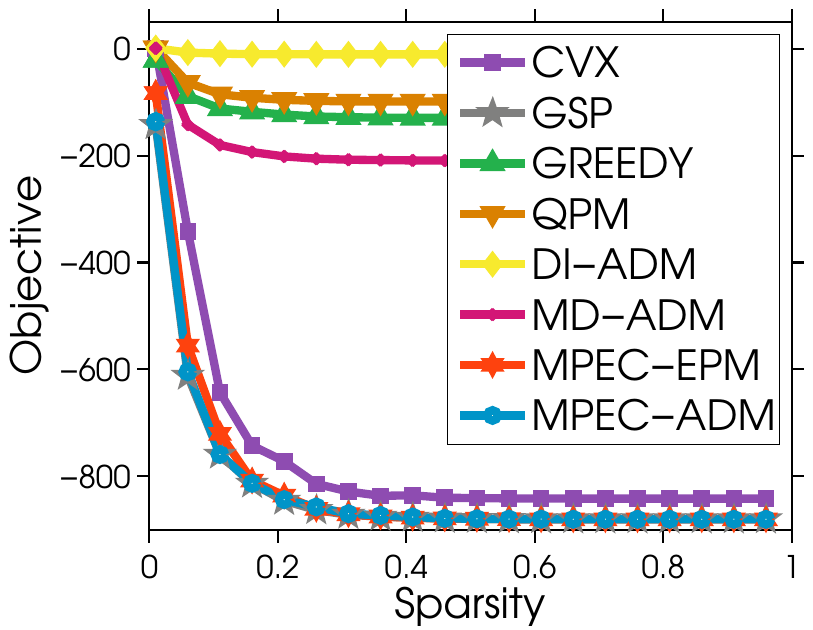}\caption{`gisette'}\end{subfigure}\ghs
      \begin{subfigure}{\imgwidlogloss}\includegraphics[width=\textwidth,height=\imgheilogloss]{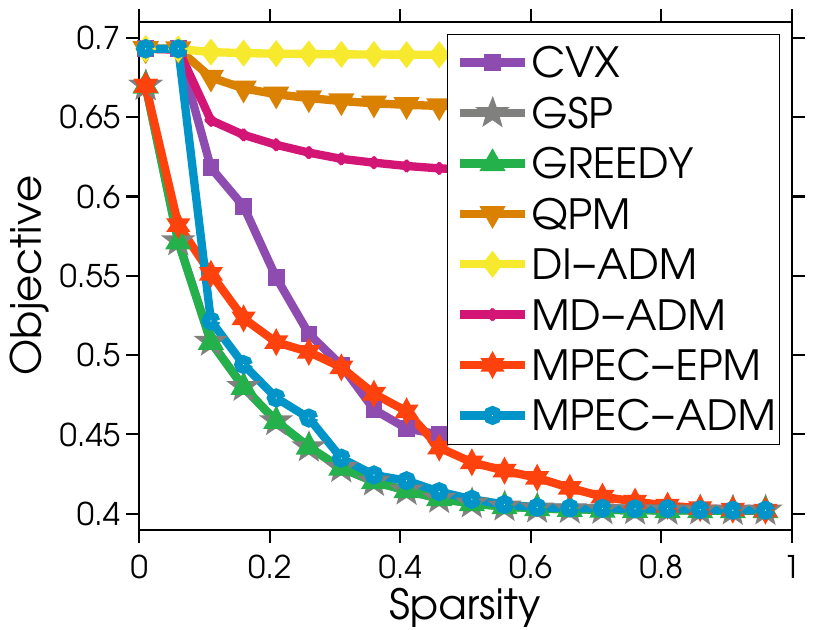}\caption{`w2a'}\end{subfigure}\ghs

      \begin{subfigure}{\imgwidlogloss}\includegraphics[width=\textwidth,height=\imgheilogloss]{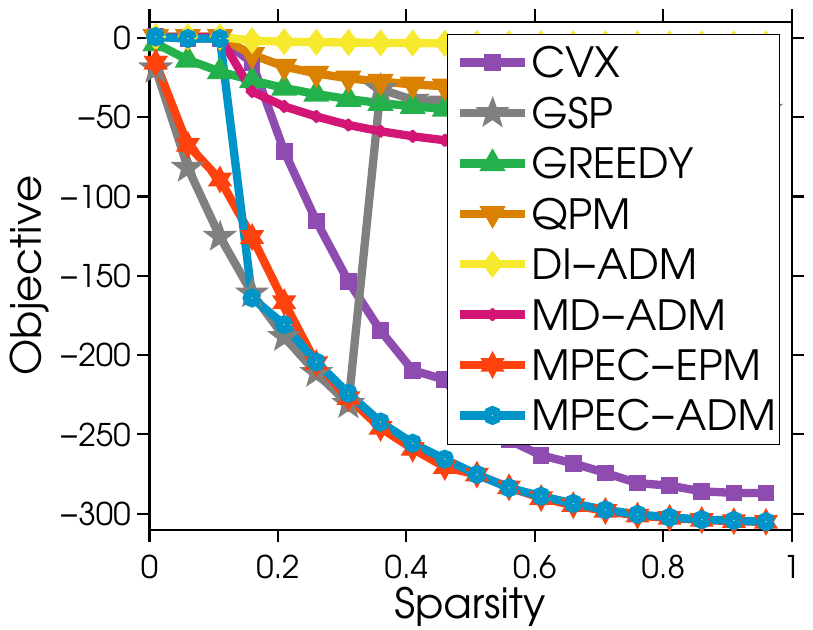}\caption{`rcv1binary'}\end{subfigure}\ghs
      \begin{subfigure}{\imgwidlogloss}\includegraphics[width=\textwidth,height=\imgheilogloss]{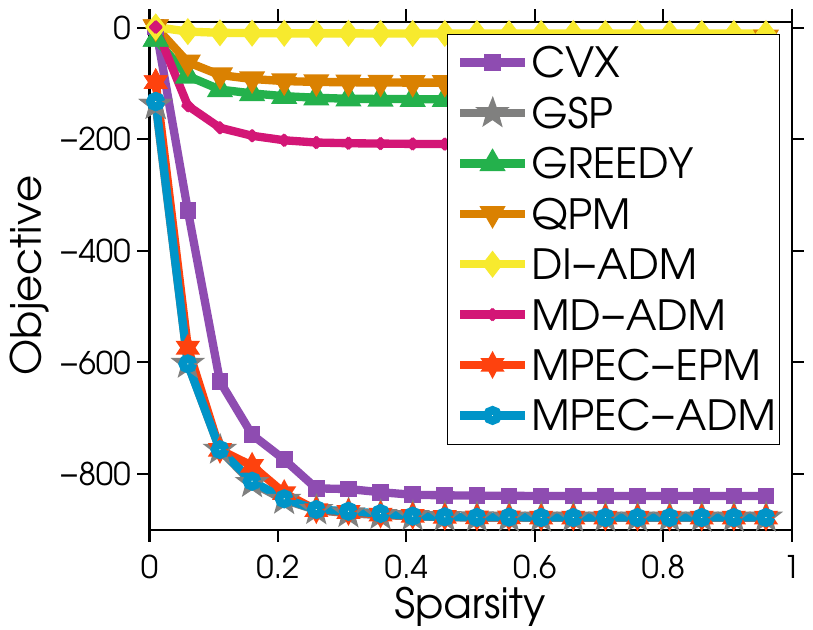}\caption{`a5a'}\end{subfigure}
      \begin{subfigure}{\imgwidlogloss}\includegraphics[width=\textwidth,height=\imgheilogloss]{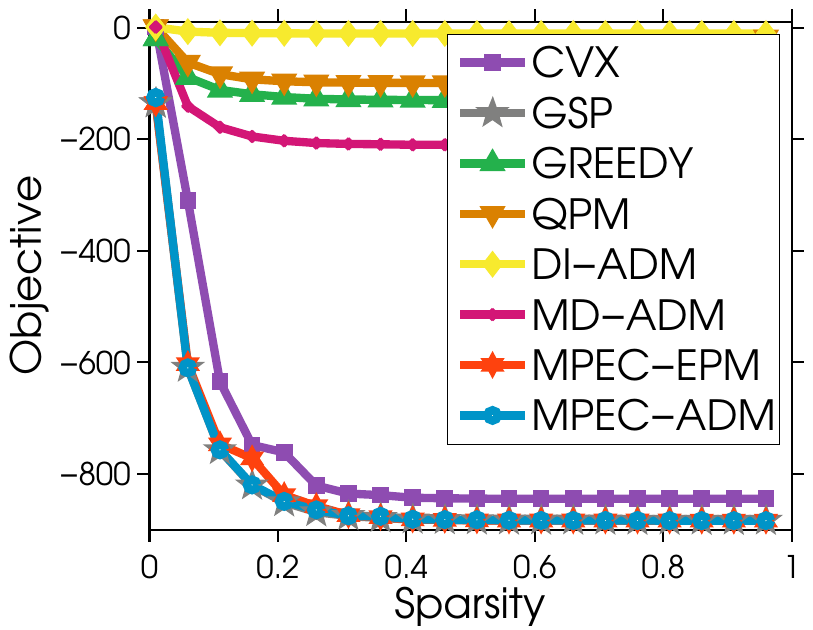}\caption{`realsim'}\end{subfigure}

      \begin{subfigure}{\imgwidlogloss}\includegraphics[width=\textwidth,height=\imgheilogloss]{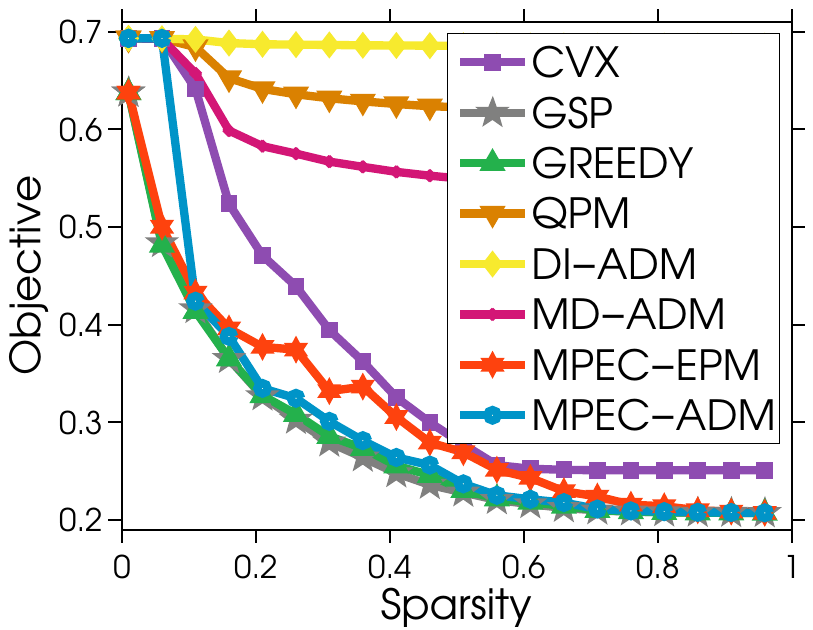}\caption{`mushrooms'}\end{subfigure}\ghs
      \begin{subfigure}{\imgwidlogloss}\includegraphics[width=\textwidth,height=\imgheilogloss]{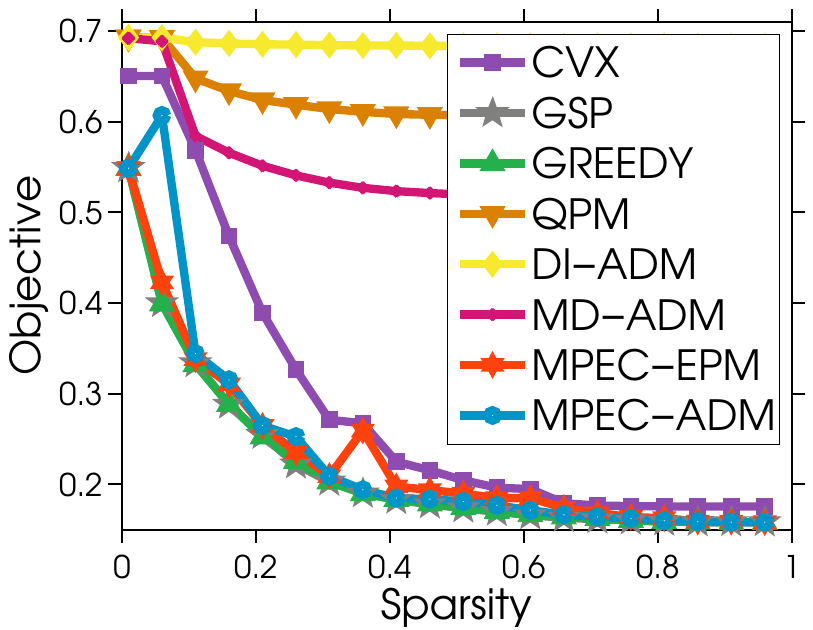}\caption{`splice'}\end{subfigure}
      \begin{subfigure}{\imgwidlogloss}\includegraphics[width=\textwidth,height=\imgheilogloss]{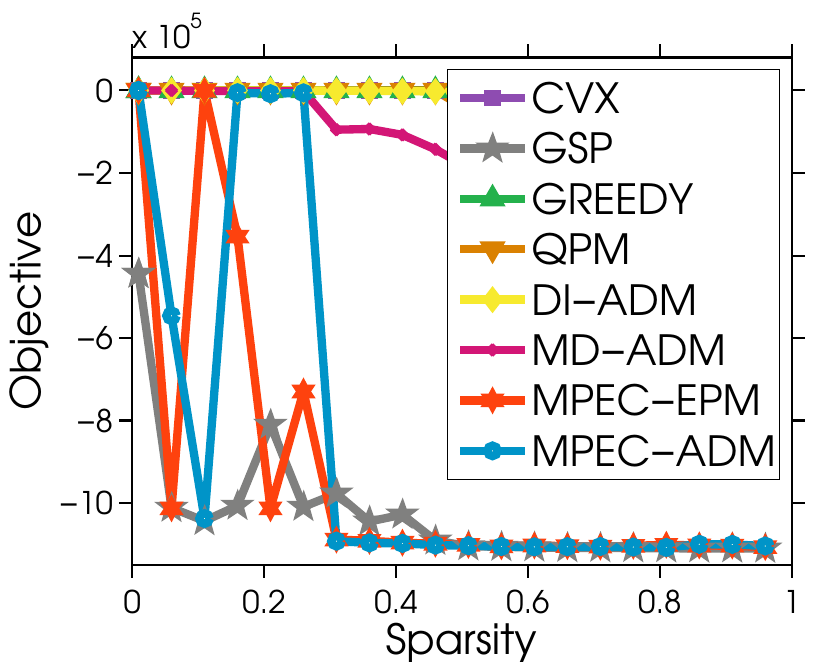}\caption{`madelon'}\end{subfigure}

      \begin{subfigure}{\imgwidlogloss}\includegraphics[width=\textwidth,height=\imgheilogloss]{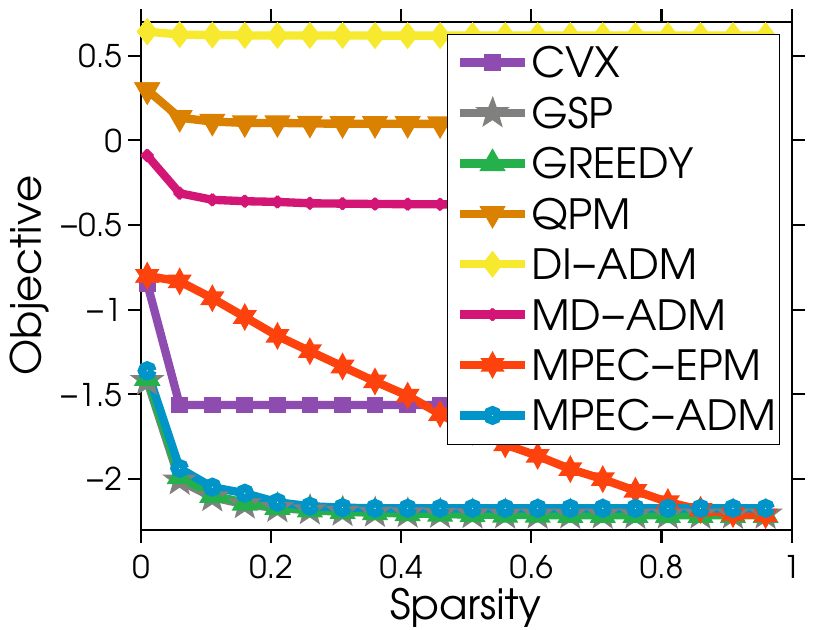}\caption{`w1a'}\end{subfigure}\ghs
      \begin{subfigure}{\imgwidlogloss}\includegraphics[width=\textwidth,height=\imgheilogloss]{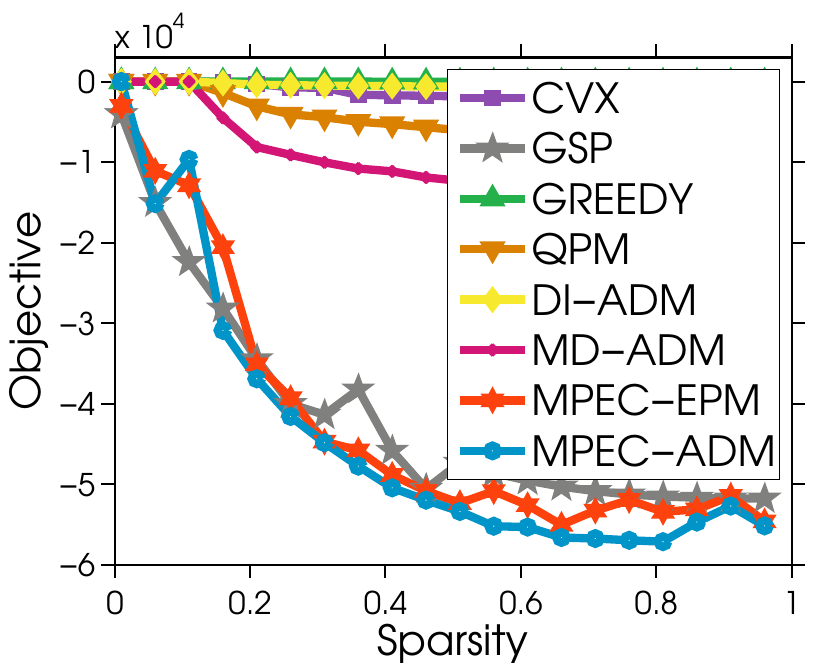}\caption{`duke'}\end{subfigure}
      \begin{subfigure}{\imgwidlogloss}\includegraphics[width=\textwidth,height=\imgheilogloss]{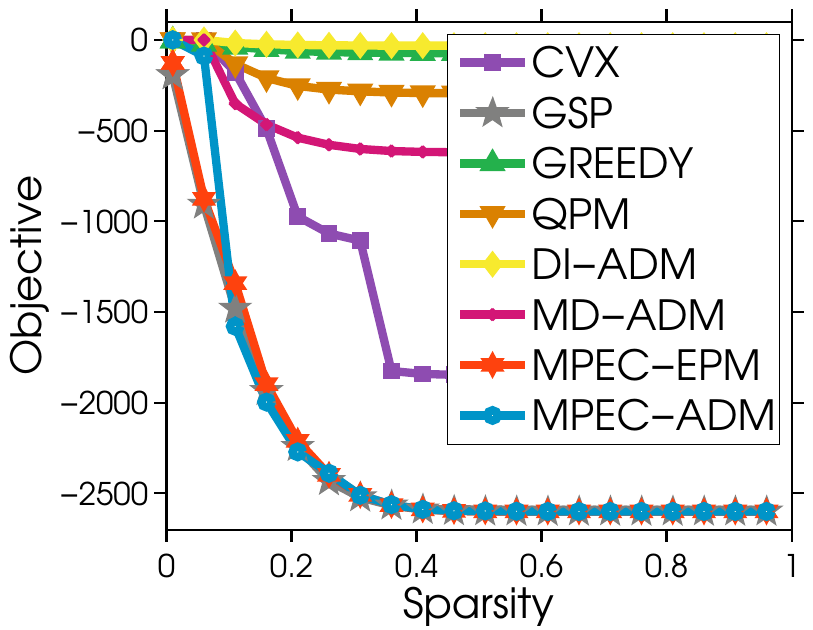}\caption{`mnist'}\end{subfigure}

\vspace{0pt}
\caption{Feature selection based on logistic loss minimization. We set $\lambda=0.01$.}\label{fig:feature:log}
\end{figure*}

\begin{figure*}[!th]
\captionsetup[subfigure]{justification=centering}
    \centering
      \begin{subfigure}{\imgwidlogloss}\includegraphics[width=\textwidth,height=\imgheilogloss]{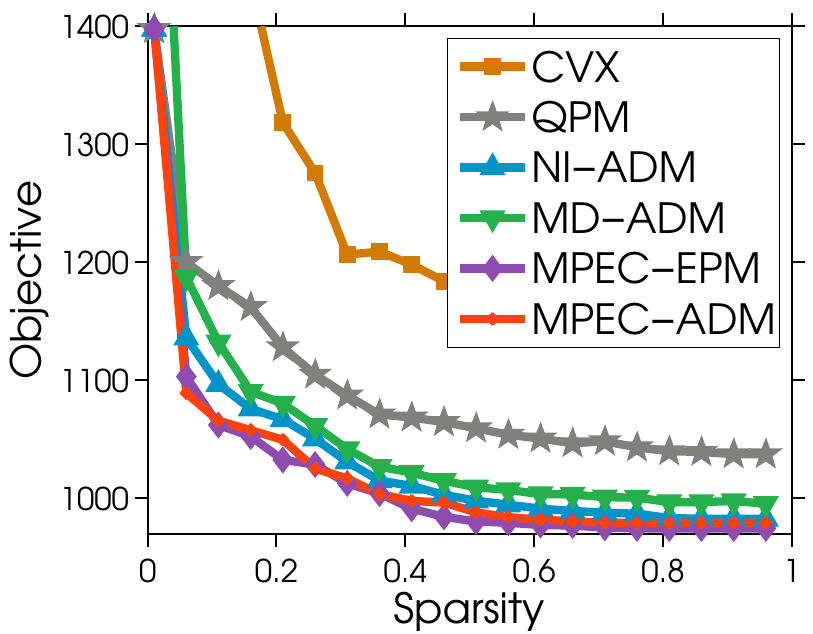}\caption{`a6a'}\end{subfigure}\ghs
      \begin{subfigure}{\imgwidlogloss}\includegraphics[width=\textwidth,height=\imgheilogloss]{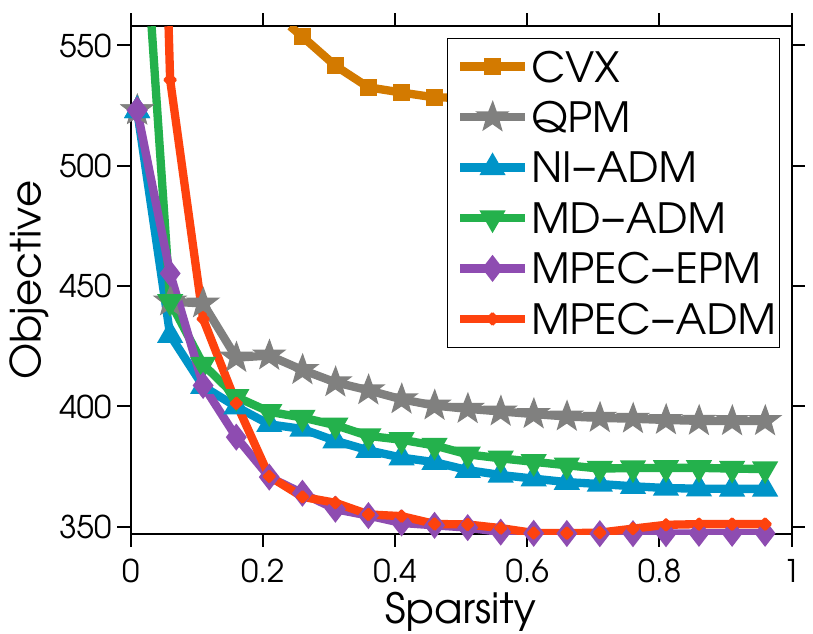}\caption{`gisette'}\end{subfigure}\ghs
      \begin{subfigure}{\imgwidlogloss}\includegraphics[width=\textwidth,height=\imgheilogloss]{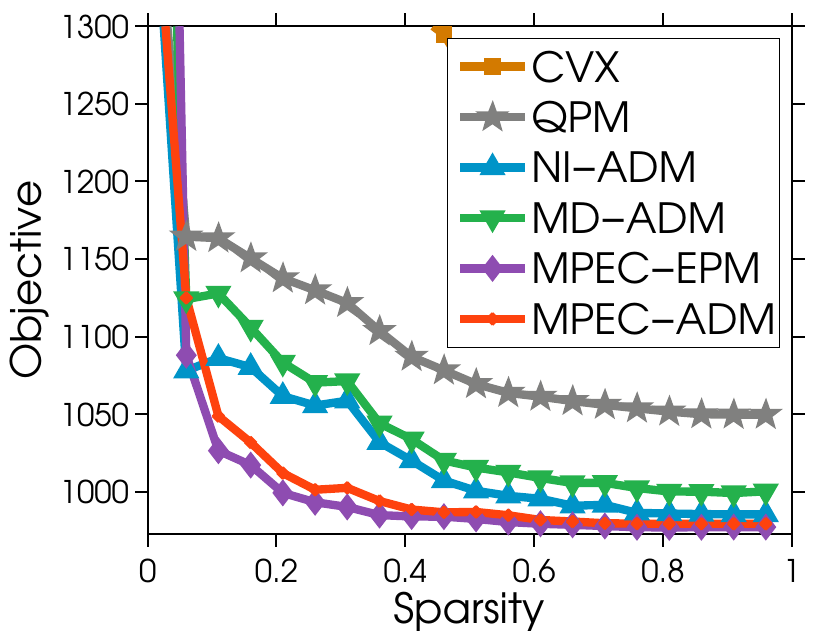}\caption{`w2a'}\end{subfigure}\ghs

      \begin{subfigure}{\imgwidlogloss}\includegraphics[width=\textwidth,height=\imgheilogloss]{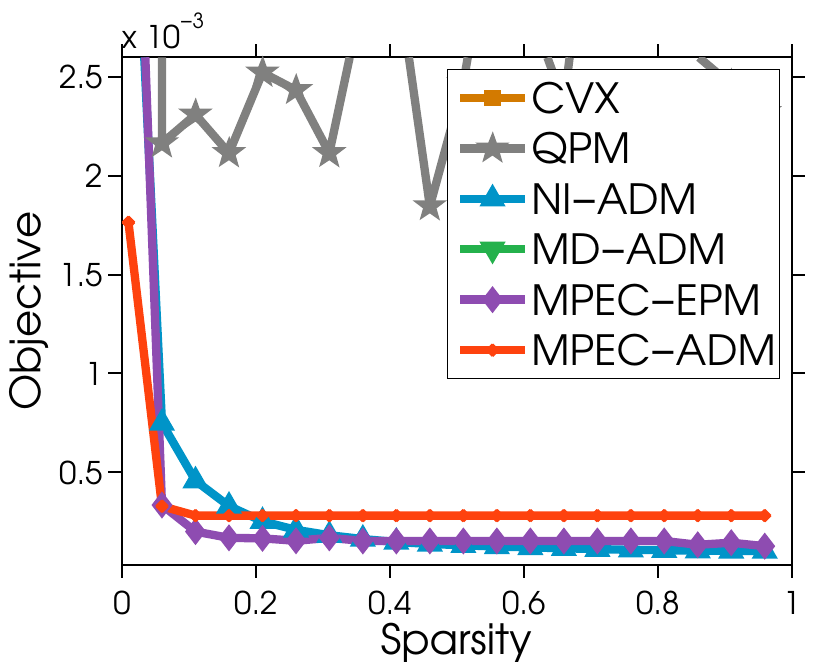}\caption{`rcv1binary'}\end{subfigure}\ghs
      \begin{subfigure}{\imgwidlogloss}\includegraphics[width=\textwidth,height=\imgheilogloss]{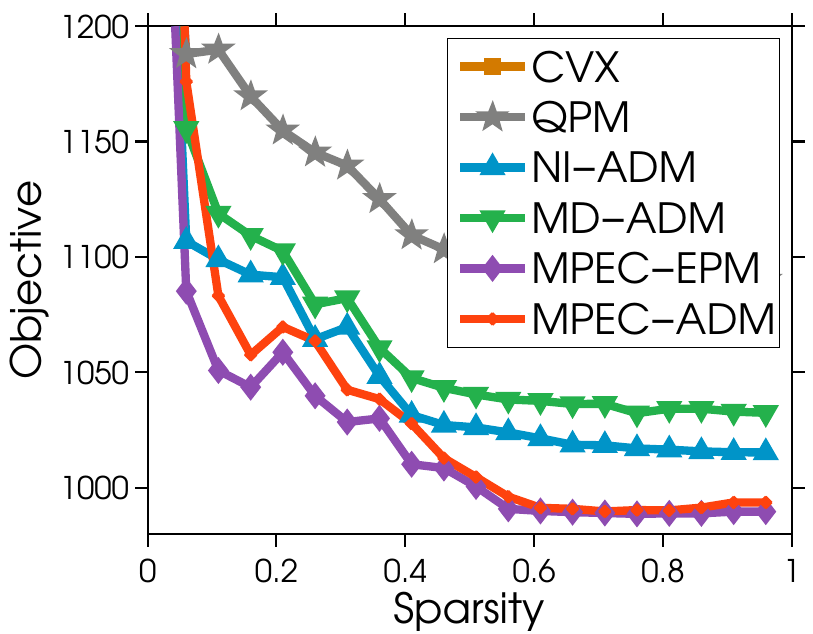}\caption{`a5a'}\end{subfigure}
      \begin{subfigure}{\imgwidlogloss}\includegraphics[width=\textwidth,height=\imgheilogloss]{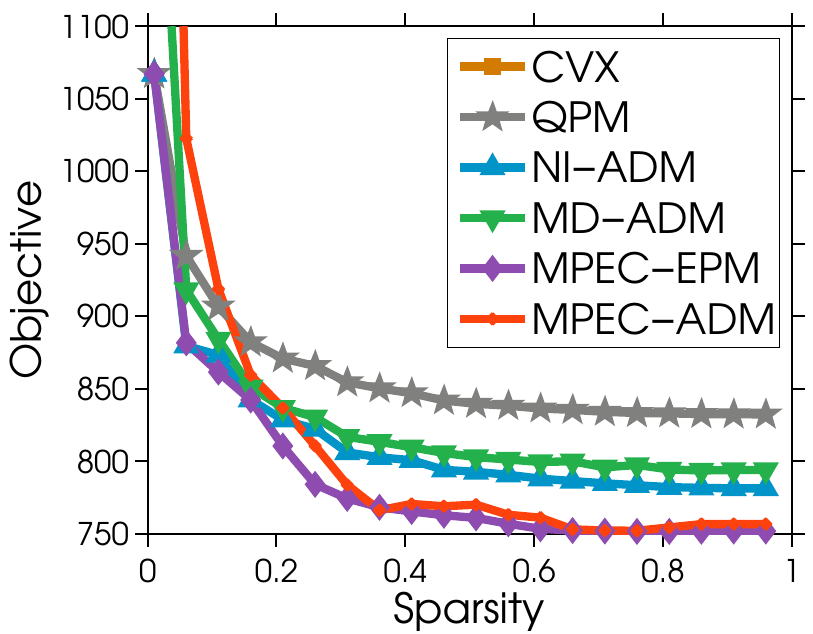}\caption{`realsim'}\end{subfigure}

      \begin{subfigure}{\imgwidlogloss}\includegraphics[width=\textwidth,height=\imgheilogloss]{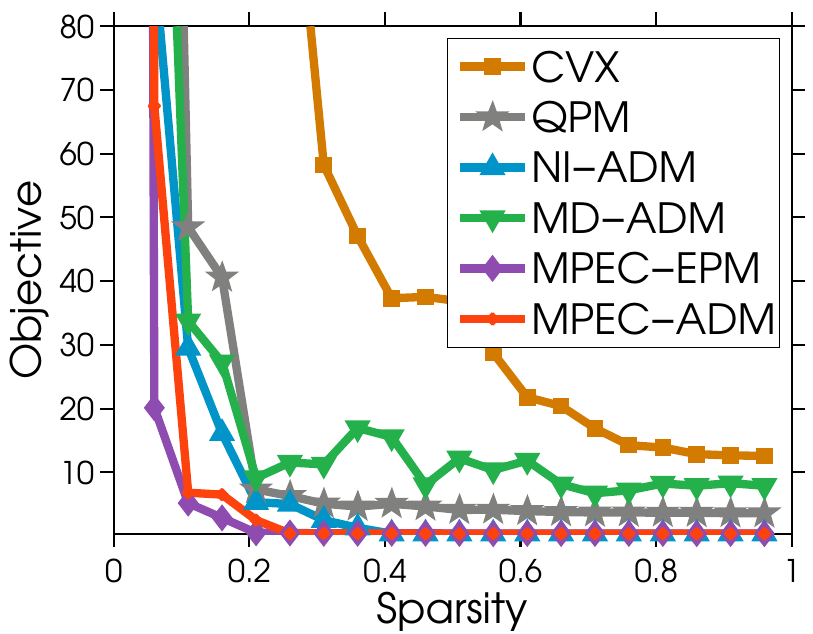}\caption{`mushrooms'}\end{subfigure}\ghs
      \begin{subfigure}{\imgwidlogloss}\includegraphics[width=\textwidth,height=\imgheilogloss]{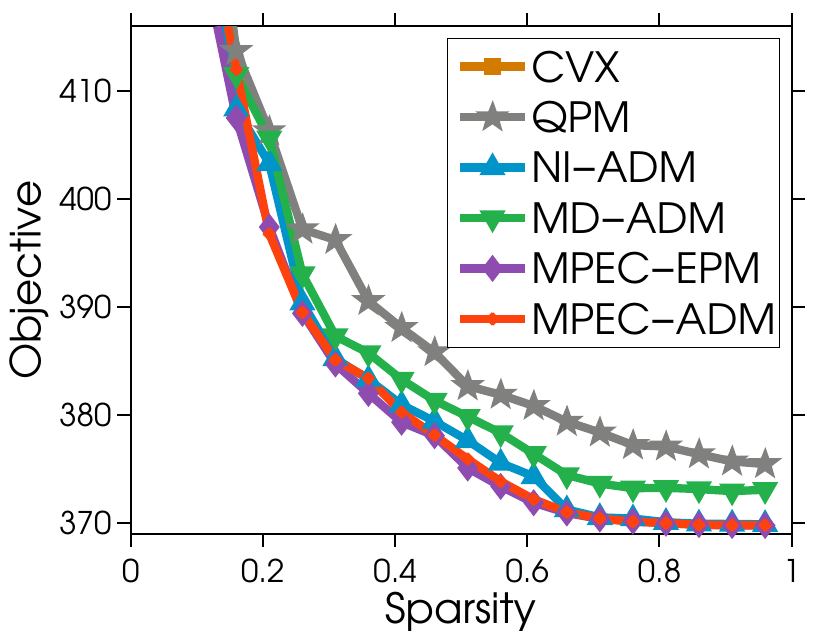}\caption{`splice'}\end{subfigure}
      \begin{subfigure}{\imgwidlogloss}\includegraphics[width=\textwidth,height=\imgheilogloss]{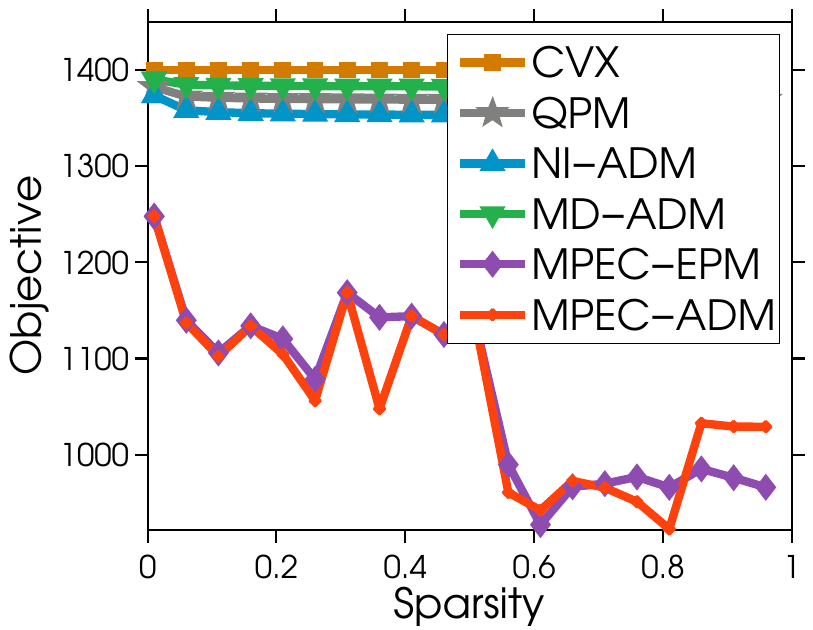}\caption{`madelon'}\end{subfigure}

      \begin{subfigure}{\imgwidlogloss}\includegraphics[width=\textwidth,height=\imgheilogloss]{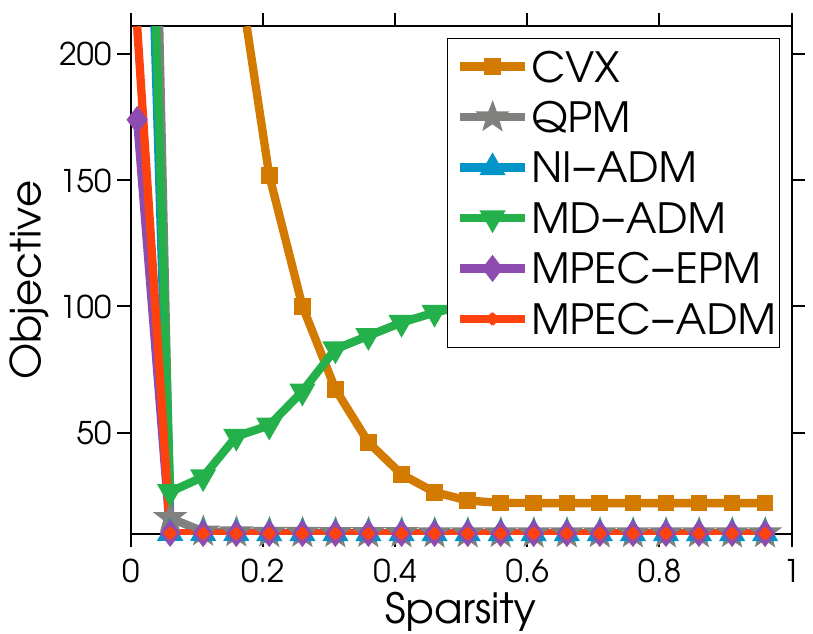}\caption{`w1a'}\end{subfigure}\ghs
      \begin{subfigure}{\imgwidlogloss}\includegraphics[width=\textwidth,height=\imgheilogloss]{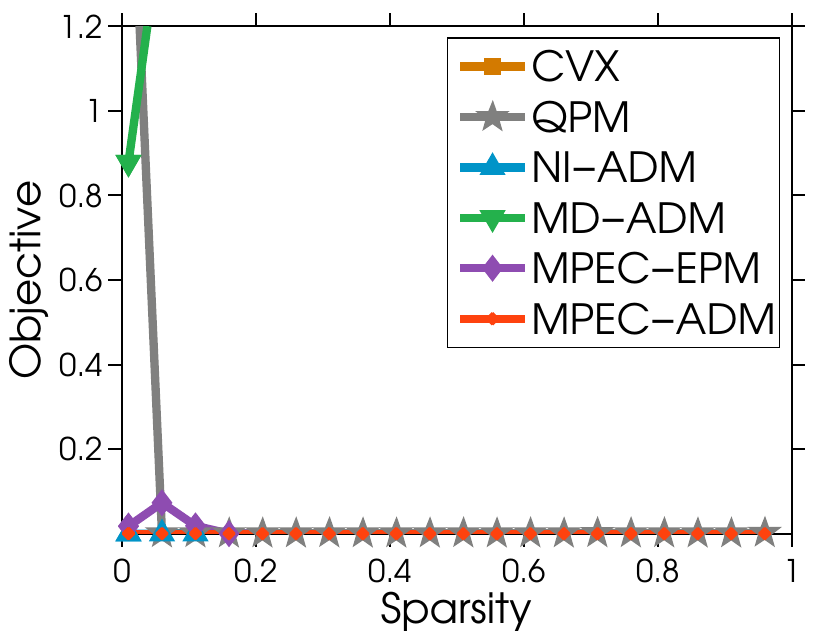}\caption{`duke'}\end{subfigure}
      \begin{subfigure}{\imgwidlogloss}\includegraphics[width=\textwidth,height=\imgheilogloss]{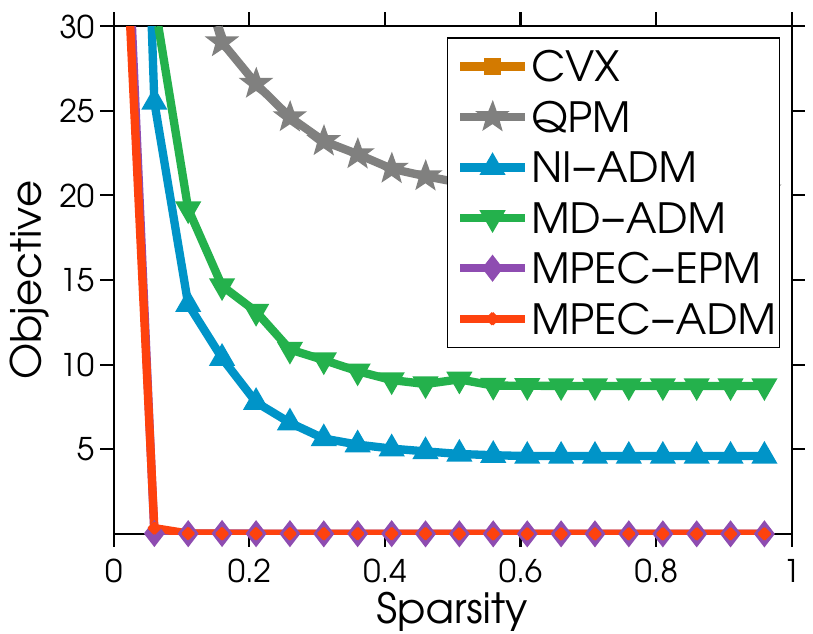}\caption{`mnist'}\end{subfigure}

\vspace{0pt}
\caption{Feature selection based on hinge loss minimization. We set $\lambda=0.01$. }
\label{fig:feature:hinge}
\end{figure*}

We include the comparisons with GREEDY and GSP only on logistic loss feature selection, since these two methods are only applicable to smooth objectives. We vary the sparsity parameter $k = \{0.01,0.06,...,0.96\} \times n$ on different data sets and report the objective value of the optimization problem\footnote{All the algorithms generally terminate within 8 minutes, we do not include the CPU time here. In fact, our methods converge within reasonable time. They are only 3-7 times slower than the classical convex methods. This is expected, since they are alternating methods and have the same computational complexity as the convex methods.}. All algorithms are implemented in Matlab on an Intel 2.6 GHz CPU with 128 GB RAM. We set $\rho^0 = 0.01$ for MPEC-EPM and $\alpha=\eta=0.01$ for MPEC-ADM. For the proximal parameter of both algorithms, we set $\mu=0.01$ \footnote{It is a small constant to guarantee strong convexity for the sub-problem and a unique solution in every iteration. In theory, the proximal strategy targets convergence in bi-linear bi-convex optimization \cite{Attouch2010,Bolte2014}. In fact, even if $\mu$ is set to zero, our algorithm is observed to converge since we apply alternating minimization on a bi-convex problem.}.


%

\subsection{Feature Selection}
Given a set of labeled patterns $\{\bbb{s}_i,\bbb{y}_i\}_{i=1}^p$, where $\bbb{s}_i$ is an instance with $n$ features and $\bbb{y}_i \in \{\pm1\}$ is the label. Feature selection considers the following optimization problem:
\beq
\min_{\bbb{x}}~\textstyle\tfrac{\lambda}{2} \|\bbb{x}\|_2^2 +  \sum_{i}^p \ell \left(   \la \bbb{s}_i,\bbb{x} \ra , \bbb{y}_i \right)  ,~s.t.~\|\bbb{x}\|_0 \leq k
\eeq
\noi where the loss function $\ell(\cdot)$ is chosen to be the logistic loss: $\ell(r,t)=\log(1+\exp(r)) - r  t$ or the hinge loss: $\ell(r,t) = \max(0,1-r t)$. 

In our experiments, we test on 12 well-known benchmark datasets\footnote{\url{www.csie.ntu.edu.tw/~cjlin/libsvmtools/datasets/}} which contain high dimensional ($n\geq 7000$) data vectors and thousands of data instances. Since we are solving an optimization problem, we only consider measuring the quality of the solutions by comparing the objective values. 

Based on our experimental results in Figure \ref{fig:feature:log}, we make the following observations. (i) Convex methods generally gives good results in this task, but they fail on `duke' and `madelon' datasets. (ii) The quadratic penalization techniques (NI-ADM, MD-ADM and QPM) generally demonstrate bad performance in this task. (iii) GREEDY and GSP often achieve good performance, but they are not stable in some cases. (iv) Both our MPEC-EPM and MPEC-ADM often achieve better performance than existing solutions. For hinge loss feature selection, we demonstrate our experimental results in Figure \ref{fig:feature:hinge}. We make the following observations. (i) Convex methods generally fail in this situation, as they usually produce much larger objectives. (ii) NI-ADM and MD-ADM outperform QPM in all test cases. (iii) MPEC-EPM and MPEC-ADM generally achieve lower objective values.

\begin{figure*}[!t]
\captionsetup[subfigure]{justification=centering}
    \centering
      \begin{subfigure}{0.3\textwidth}\includegraphics[width=\textwidth,height=\imgheiseg]{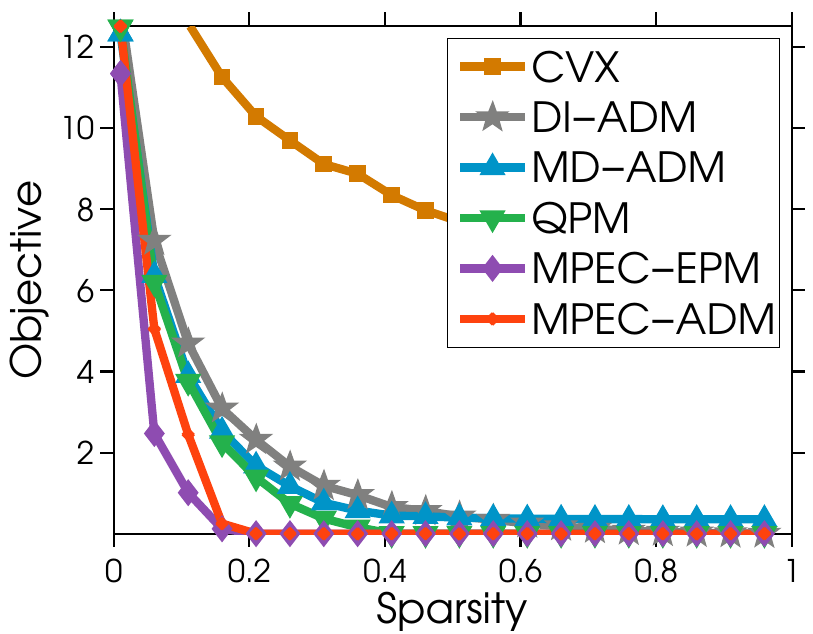}\caption{\figlabsize $n=1024$,~$\sigma=10$}\end{subfigure}\ghs
      \begin{subfigure}{0.3\textwidth}\includegraphics[width=\textwidth,height=\imgheiseg]{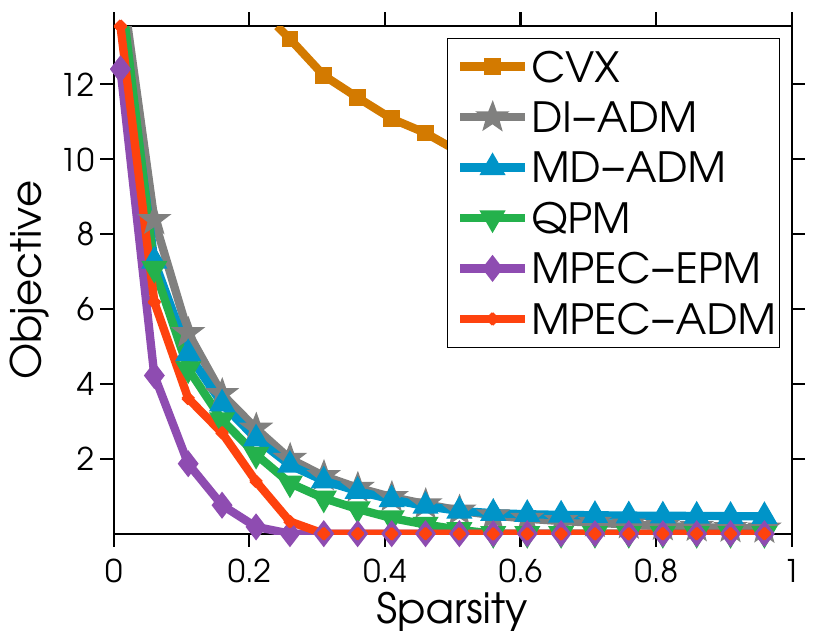}\caption{\figlabsize $n=2048$,~$\sigma=10$}\end{subfigure}\ghs
      \begin{subfigure}{0.3\textwidth}\includegraphics[width=\textwidth,height=\imgheiseg]{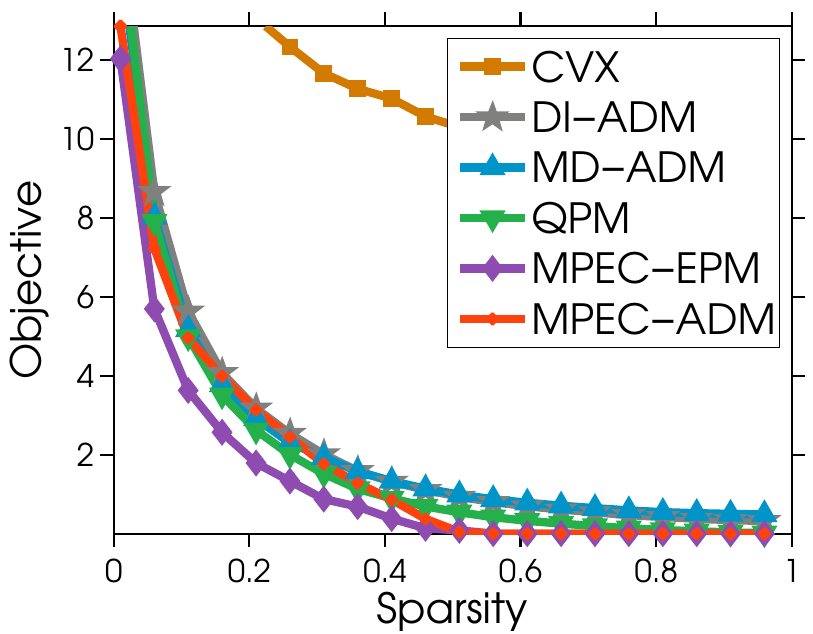}\caption{\figlabsize $n=4096$,~$\sigma=10$}\end{subfigure}

      \begin{subfigure}{0.3\textwidth}\includegraphics[width=\textwidth,height=\imgheiseg]{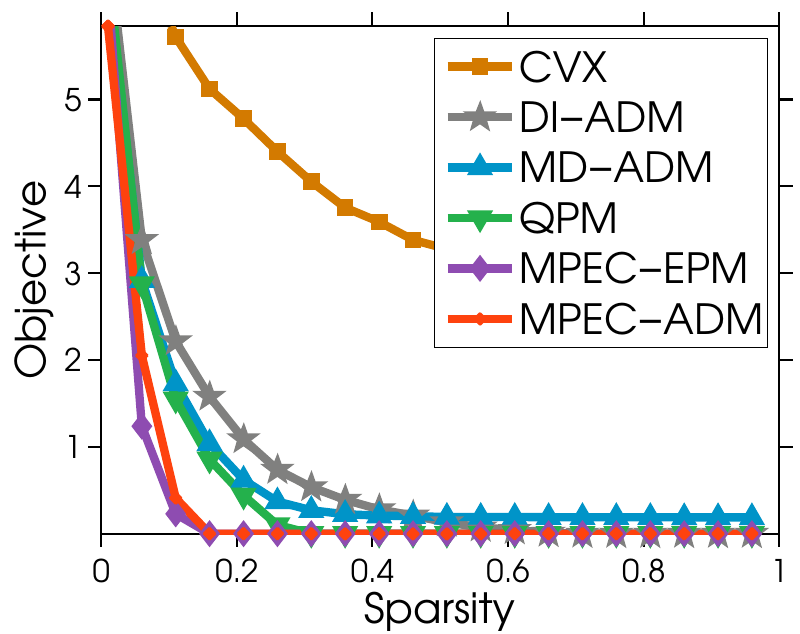}\caption{\figlabsize $n=1024$,~$\sigma=5$}\end{subfigure}\ghs
      \begin{subfigure}{0.3\textwidth}\includegraphics[width=\textwidth,height=\imgheiseg]{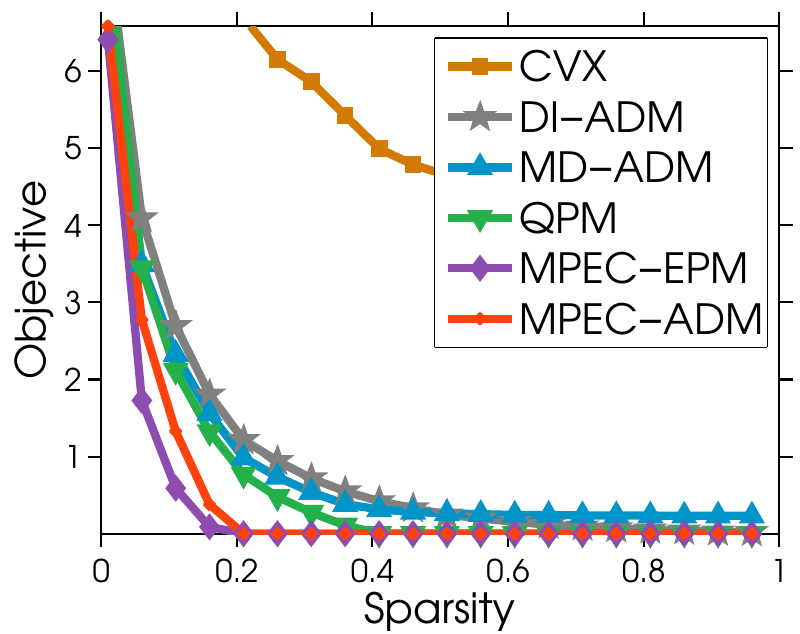}\caption{\figlabsize $n=2048$,~$\sigma=5$}\end{subfigure}\ghs
      \begin{subfigure}{0.3\textwidth}\includegraphics[width=\textwidth,height=\imgheiseg]{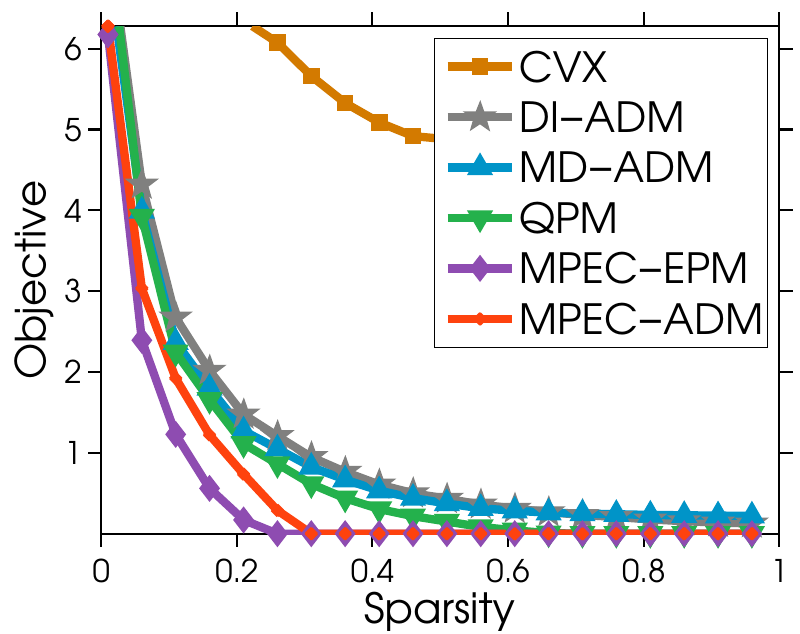}\caption{\figlabsize $n=4096$,~$\sigma=5$}\end{subfigure}
\vspace{0pt}
\caption{Segmented regression. Results for unit-column-norm $\bbb{A}$ for different parameter settings.}
\label{fig:dan:sel}
\end{figure*}

\begin{figure*} [!t]
\begin{center}
\includegraphics[width=0.8\textwidth,height=8cm]{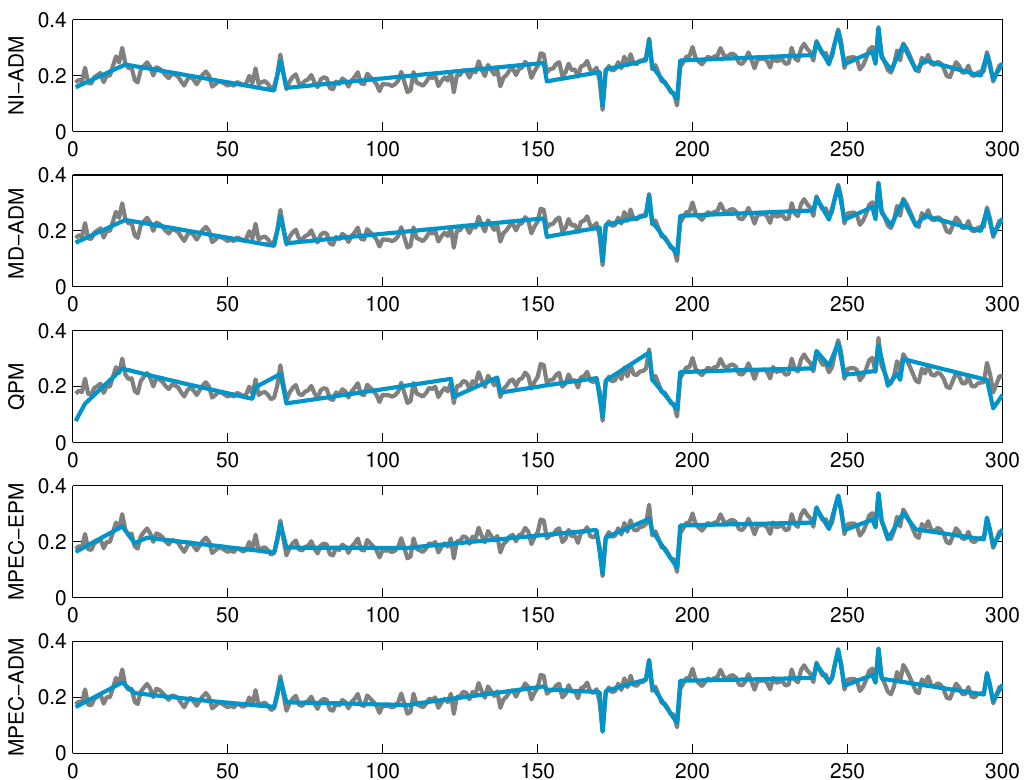}
\caption{Trend filtering result on `snp500.dat' data \cite{KimKBG09}. The objective values achieved by the five methods are $0.38, 0.38, 0.53, 0.32, 0.33$, respectively.}
\label{fig:trend}
\end{center}
\end{figure*}

\subsection{Segmented Regression}

Given a design matrix $\bbb{A}\in \mathbb{R}^{m\times n}$ and an observation vector $\bbb{b}\in \mathbb{R}^{m}$, segmented regression involves solving the following optimization problem:
\beq \label{eq:seg:reg}
\min_{\bbb{x}}~\textstyle\| \bbb{A}^T(\bbb{A}x-\bbb{b}) \|_{\infty},~s.t.~\|\bbb{x}\|_0\leq k
\eeq
\noi This regression model is closely related to Dantzig selector \cite{candes2007dantzig} in the literature, where $\ell_0$ norm is replaced by $\ell_1$ norm in Eq (\ref{eq:seg:reg}). In our experiments, we consider design matrices with unit column norms. Similar to \cite{candes2007dantzig}, we first generate an $\bbb{A}$ with independent Gaussian entries and then normalize each column to have norm 1. We then select a support set of size $0.5 \times \min(m,n)$ uniformly at random. We finally set $\bbb{b} = \bbb{Ax} + \bbb{\varepsilon} $ with $\bbb{\varepsilon} \sim N(0,\sigma^2 \bbb{I})$. We vary $n$ from $\{1024,2048,4096\}$ and set $m = n/8$. We demonstrate our experimental results in Figure \ref{fig:dan:sel}. The proposed penalization techniques (MPEC-EPM and MPEC-ADM) outperform the quadratic penalization techniques (NI-ADM, MD-ADM and QPM) consistently.

\subsection{Trend Filtering}\label{exp:tf}
Given a time series $\bbb{y}\in \mathbb{R}^n$, the goal of trend filtering is to find another nearest time series $\bbb{x}$ such that $\bbb{x}$ is sparse after a gradient mapping \cite{KimKBG09}. Trend filtering involves solving the following optimization, where $\bbb{D} \in \mathbb{R}^{(n-2)\times n}$ is a $2^{\text{nd}}$ difference matrix.
\beq
 \min_{\bbb{x}}~\textstyle\tfrac{1}{2}\|\bbb{x}-\bbb{y}\|_2^2, ~s.t.~\|\bbb{Dx}\|_0 \leq k
\eeq

In our experiments, we perform trend filtering on the `snp500.dat' data set \cite{KimKBG09}. For better visualization, we only use the first $300$ time instances of the signal $\bbb{y}$, i.e. $n=300$. The sparsity parameter $k$ is set to $30$. We stop all the algorithms when the same strict stopping criterion is met, in order to ensure that the hard constraint $\|\bbb{Ax}\|_0 \leq k$ is fully satisfied. As shown in Figure \ref{fig:trend}, the proposed methods (MPEC-EPM and MPEC-ADM) achieve more natural trend filtering results (refer to position 120-170), since they achieve lower objective values (0.32 and 0.33).

\begin{figure*}[!t]
\captionsetup[subfigure]{justification=centering}
    \centering
      \begin{subfigure}{\imgwidmrf}\includegraphics[width=\textwidth,height=\imgheimrf]{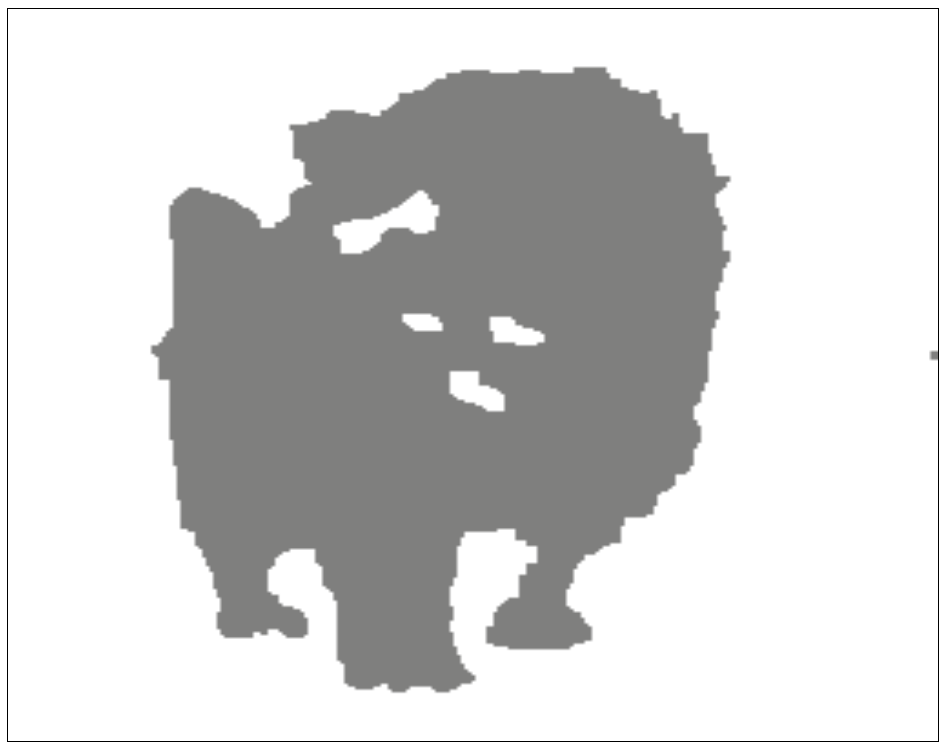}\caption{\figlabsize by Graph Cut\\ $f(\mathbf{x})=95741.1$}\end{subfigure}\ghs
      \begin{subfigure}{\imgwidmrf}\includegraphics[width=\textwidth,height=\imgheimrf]{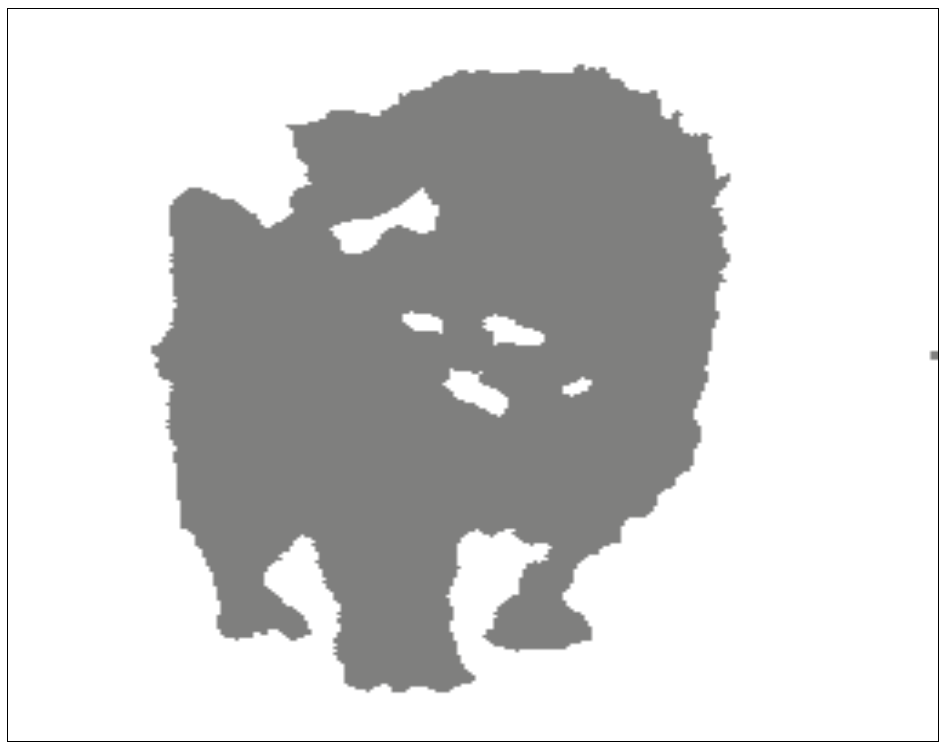}\caption{\figlabsize by NI-ADM\\ $f(\mathbf{x})=95852.6$}\end{subfigure}\ghs
      \begin{subfigure}{\imgwidmrf}\includegraphics[width=\textwidth,height=\imgheimrf]{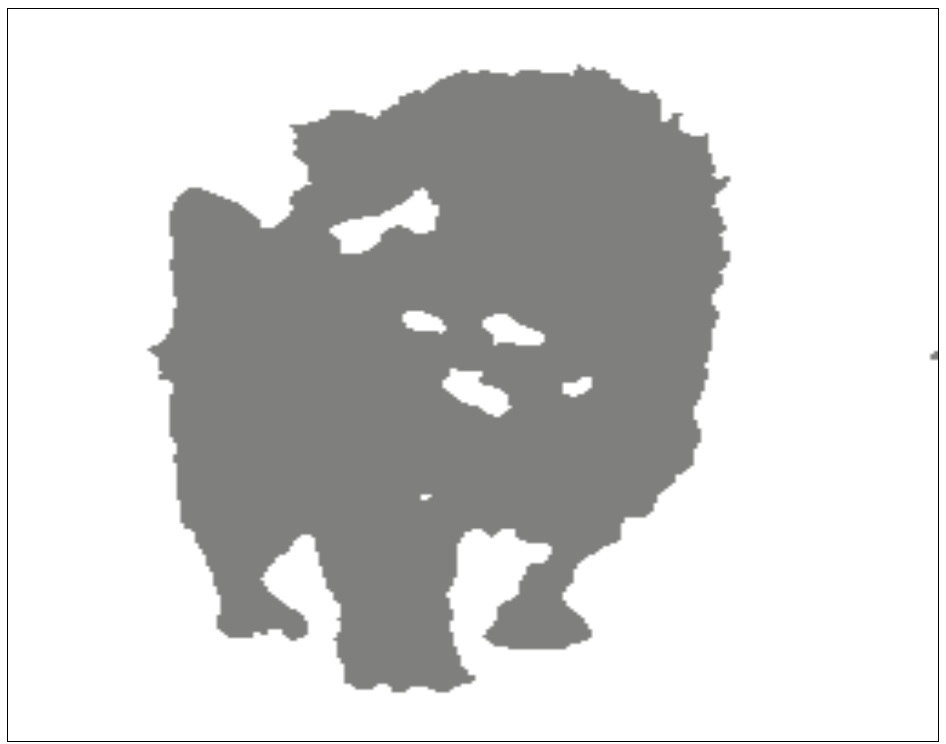}\caption{\figlabsize by MD-ADM\\ $f(\mathbf{x})=95825.2$}\end{subfigure}

      \begin{subfigure}{\imgwidmrf}\includegraphics[width=\textwidth,height=\imgheimrf]{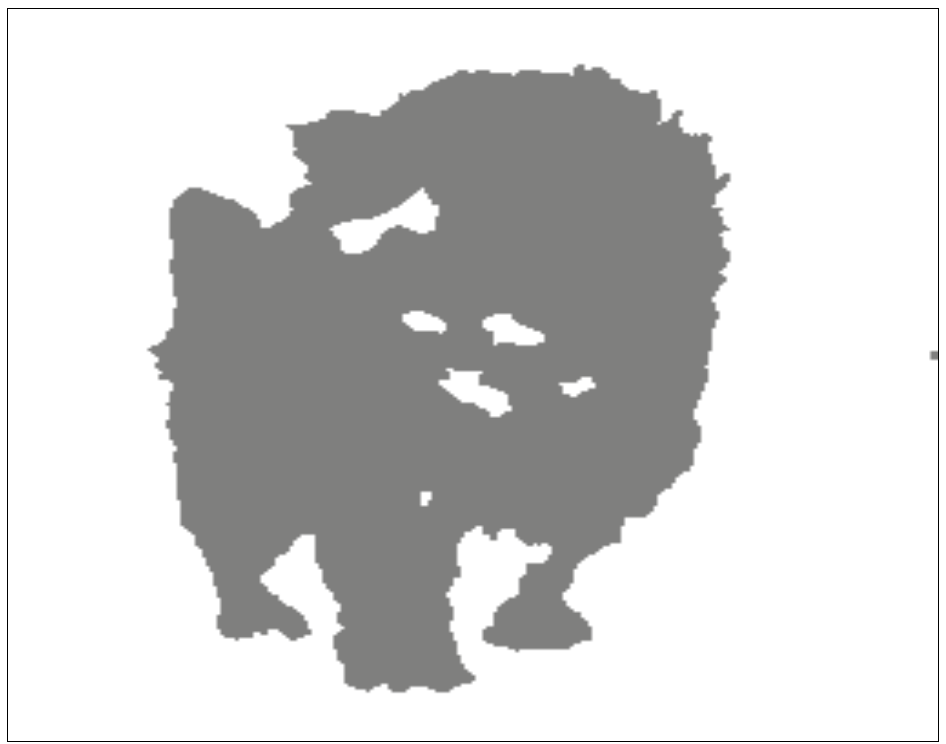}\caption{\figlabsize by QPM\\ $f(\mathbf{x})=95827.3$}\end{subfigure}\ghs
      \begin{subfigure}{\imgwidmrf}\includegraphics[width=\textwidth,height=\imgheimrf]{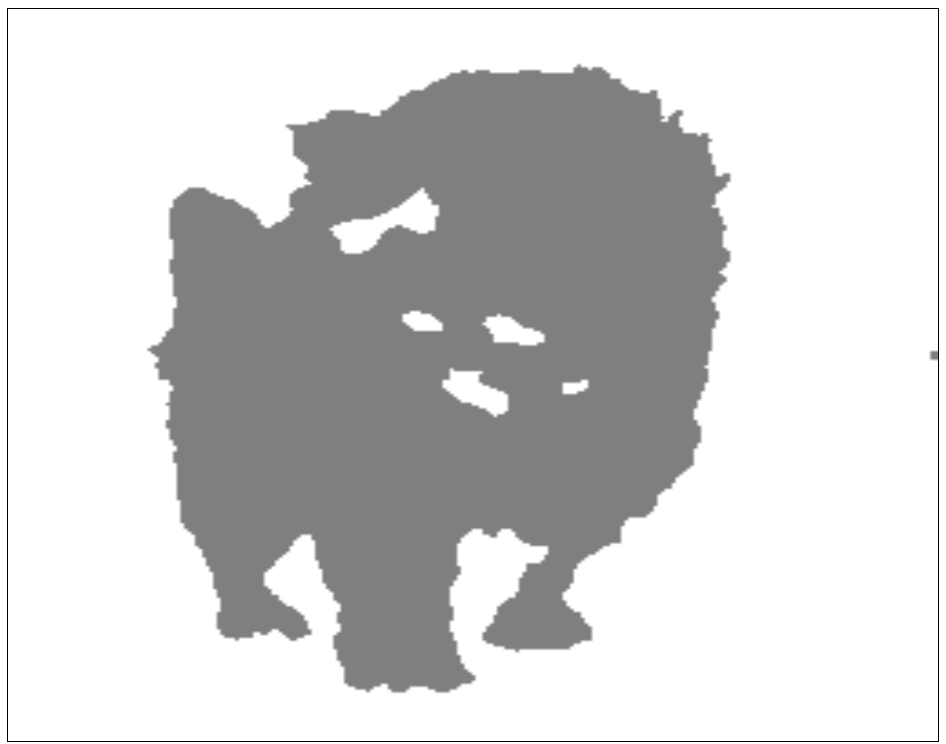}\caption{\figlabsize by MPEC-EPM\\ $f(\mathbf{x})=95798.8$}\end{subfigure}\ghs
      \begin{subfigure}{\imgwidmrf}\includegraphics[width=\textwidth,height=\imgheimrf]{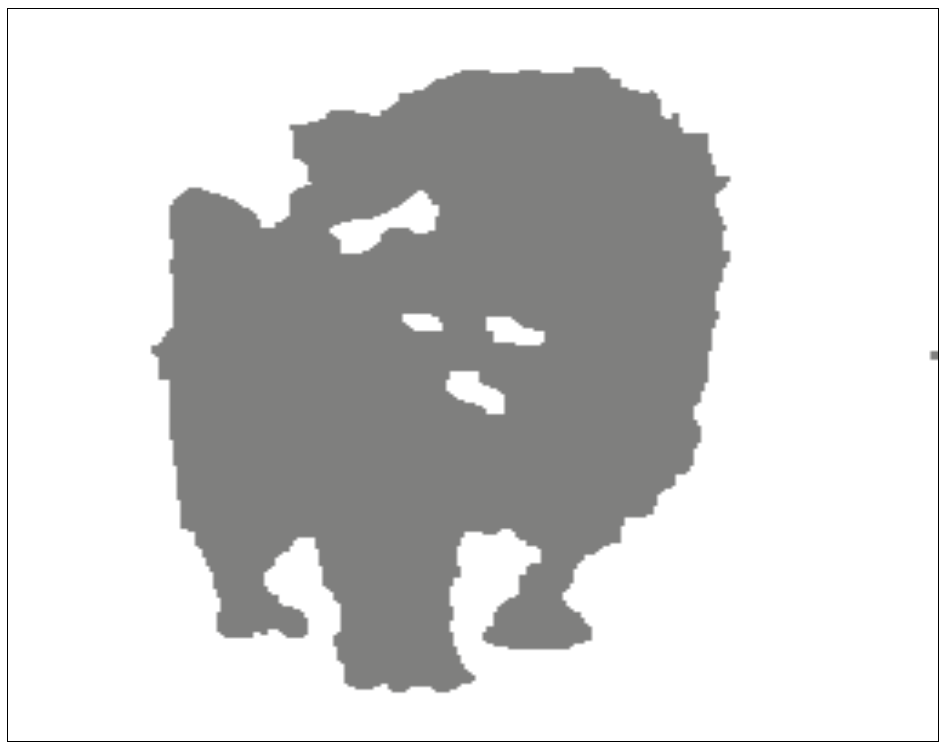}\caption{\figlabsize by MPEC-ADM\\ $f(\mathbf{x})=95748.0$}\end{subfigure}

\caption{MRF optimization on `cat' image. The objective value for the unary image $\bbb{b}$ is $f(\bbb{b})=96439.8$. }
\label{fig:mrf}
\end{figure*}

\subsection{MRF Optimization} \label{exp:mrf}
Markov Random Field (MRF) optimization involves solving the following problem \cite{KolmogorovZ04},
\beq
\min_{\bbb{x}}~\textstyle \frac{1}{2} \bbb{x}^T\bbb{Lx}+\bbb{b}^T\bbb{x},~s.t.~\bbb{x}\in\{0,1\}^n
\eeq
\noi where $\bbb{b}\in\mathbb{R}^n$ is determined by the unary term defined for the graph and $\bbb{L}\in\mathbb{R}^{n\times n}$ is the Laplacian matrix, which is based on the binary term relating pairs of graph nodes together. The quadratic term is usually considered a smoothness prior on the node labels. The MRF formulation is widely used in many labeling applications, including image segmentation.

In our experiments, we include the comparison with the graph cut method \cite{KolmogorovZ04}, which is known to achieve the global optimal solution for this specific class of binary problem. Figure \ref{fig:mrf} demonstrates a qualitative result for image segmentation. MPEC-ADM produces a solution that is very close to the globally optimal one. Moreover, both our methods achieve lower objectives than the other compared methods.

\begin{table*}[!t]
\tiny\small \footnotesize
\begin{center}
\caption{Random-value impulse noise removal problems. The results separated by `/' are $SNR_{0}$, $SNR_{1}$ and $SNR_{2}$, respectively. }
\label{denoising:rv}
\scalebox{0.82}{\begin{tabular}{|p{1.8cm}|p{2.2cm}|p{2.2cm}|p{2.2cm}|p{2.2cm}|p{2.2cm}|p{2.2cm}|}
\hline
\diagbox{Img.}{Alg.} &  CVX (i.e. $\ell_1TV$) & DI-ADM & MD-ADM            & QPM & MPEC-EPM & MPEC-ADM                \\
\hline
\hline
mandrill+10\%  &  0.83/3.06/5.00  &  0.78/3.68/7.06  &  0.84/3.71/7.10  &  0.85/3.76/7.16  &  0.85/3.71/6.90  &  \textcolor{red}{\textbf{0.86}}/\textcolor{red}{\textbf{3.87}}/\textcolor{red}{\textbf{7.64}}  \\
mandrill+30\%  &  0.69/2.38/4.42  &  0.67/2.96/\textcolor{red}{\textbf{5.13}}  &  0.74/2.91/4.96  &  0.75/2.94/5.02  &  0.75/2.89/4.87  &  \textcolor{red}{\textbf{0.76}}/\textcolor{red}{\textbf{3.02}}/\textcolor{red}{\textbf{5.13}}  \\
mandrill+50\%  &  0.55/1.80/3.33  &  0.58/2.32/3.86  &  0.65/2.23/3.64  &  0.66/2.30/3.79  &  0.66/2.22/3.58  &  \textcolor{red}{\textbf{0.67}}/\textcolor{red}{\textbf{2.39}}/\textcolor{red}{\textbf{3.88}}  \\
mandrill+70\%  &  0.41/0.91/1.64  &  0.47/1.49/2.42  &  0.52/1.39/2.17  &  0.52/1.40/2.22  &  0.53/1.39/2.12  &  \textcolor{red}{\textbf{0.55}}/\textcolor{red}{\textbf{1.56}}/\textcolor{red}{\textbf{2.43}}  \\
mandrill+90\%  &  0.31/0.02/0.09  &  0.33/0.14/\textcolor{red}{\textbf{0.13}}  &  0.35/0.04/-0.08  &  0.35/0.04/-0.03  &  0.36/0.07/-0.14  &  \textcolor{red}{\textbf{0.37}}/\textcolor{red}{\textbf{0.19}}/0.05  \\
\hline
lenna+10\%  &  0.97/8.08/13.48  &  0.97/8.05/13.51  &  0.97/8.13/13.59  &  0.97/8.12/13.87  &  0.97/8.10/13.72  &  \textcolor{red}{\textbf{0.98}}/\textcolor{red}{\textbf{8.31}}/\textcolor{red}{\textbf{14.55}}  \\
lenna+30\%  &  0.92/7.40/13.60  &  0.92/6.71/10.53  &  0.91/6.62/10.07  &  0.94/7.38/12.27  &  0.93/7.19/11.86  &  \textcolor{red}{\textbf{0.96}}/\textcolor{red}{\textbf{7.90}}/\textcolor{red}{\textbf{13.68}}  \\
lenna+50\%  &  0.84/5.50/9.46  &  0.84/5.32/8.13  &  0.83/5.18/7.53  &  0.88/6.28/10.10  &  0.86/5.75/8.69  &  \textcolor{red}{\textbf{0.91}}/\textcolor{red}{\textbf{6.65}}/\textcolor{red}{\textbf{10.89}}  \\
lenna+70\%  &  0.57/2.48/4.17  &  0.67/3.34/5.11  &  0.66/3.17/4.40  &  0.71/3.66/5.34  &  0.69/3.34/4.53  &  \textcolor{red}{\textbf{0.75}}/\textcolor{red}{\textbf{4.09}}/\textcolor{red}{\textbf{5.81}}  \\
lenna+90\%  &  0.35/0.53/0.89  &  0.41/0.76/\textcolor{red}{\textbf{0.91}}  &  0.40/0.54/0.48  &  0.41/0.66/0.80  &  0.43/0.54/0.38  &  \textcolor{red}{\textbf{0.48}}/\textcolor{red}{\textbf{0.89}}/0.80  \\
\hline
lake+10\%  &  0.94/8.03/14.66  &  0.95/8.21/13.64  &  \textcolor{red}{\textbf{0.96}}/8.31/13.81  &  \textcolor{red}{\textbf{0.96}}/8.40/14.52  &  \textcolor{red}{\textbf{0.96}}/8.33/14.14  &  \textcolor{red}{\textbf{0.96}}/\textcolor{red}{\textbf{8.45}}/\textcolor{red}{\textbf{14.72}}  \\
lake+30\%  &  0.88/7.29/12.83  &  0.86/6.86/10.87  &  0.87/6.81/10.02  &  0.91/7.57/12.45  &  0.89/7.26/11.33  &  \textcolor{red}{\textbf{0.92}}/\textcolor{red}{\textbf{7.87}}/\textcolor{red}{\textbf{13.34}}  \\
lake+50\%  &  0.61/4.89/8.31  &  0.76/5.69/8.52  &  0.77/5.65/8.02  &  0.83/6.57/10.36  &  0.80/5.98/8.78  &  \textcolor{red}{\textbf{0.85}}/\textcolor{red}{\textbf{6.75}}/\textcolor{red}{\textbf{10.56}}  \\
lake+70\%  &  0.40/2.07/3.61  &  0.51/3.68/5.49  &  0.56/3.72/5.14  &  0.60/4.01/5.60  &  0.61/3.94/5.28  &  \textcolor{red}{\textbf{0.62}}/\textcolor{red}{\textbf{4.08}}/\textcolor{red}{\textbf{5.61}}  \\
lake+90\%  &  0.27/0.54/0.87  &  0.25/\textcolor{red}{\textbf{0.94}}/\textcolor{red}{\textbf{1.06}}  &  0.27/0.87/0.65  &  0.26/0.75/0.65  &  \textcolor{red}{\textbf{0.29}}/0.81/0.44  &  0.28/0.90/0.71  \\
\hline
jetplane+10\%  &  \textcolor{red}{\textbf{0.42}}/\textcolor{red}{\textbf{3.33}}/\textcolor{red}{\textbf{7.98}}  &  0.37/3.16/7.10  &  0.39/3.16/7.15  &  0.39/3.20/7.43  &  0.39/3.16/7.28  &  0.39/3.22/7.52  \\
jetplane+30\%  &  0.36/2.94/7.00  &  0.35/2.70/5.69  &  0.38/2.68/5.60  &  0.39/2.98/6.68  &  0.39/2.87/6.24  &  \textcolor{red}{\textbf{0.40}}/\textcolor{red}{\textbf{3.14}}/\textcolor{red}{\textbf{7.23}}  \\
jetplane+50\%  &  0.26/1.75/4.24  &  0.31/2.02/4.21  &  0.37/2.08/4.00  &  0.38/2.56/5.45  &  0.37/2.26/4.48  &  \textcolor{red}{\textbf{0.40}}/\textcolor{red}{\textbf{2.72}}/\textcolor{red}{\textbf{5.75}}  \\
jetplane+70\%  &  0.24/-0.55/-0.17  &  0.24/0.62/1.46  &  0.31/0.83/1.49  &  0.30/0.81/1.44  &  \textcolor{red}{\textbf{0.34}}/\textcolor{red}{\textbf{0.98}}/\textcolor{red}{\textbf{1.55}}  &  0.33/0.92/1.38  \\
jetplane+90\%  &  0.19/-1.81/-3.38  &  0.17/-1.98/-3.23  &  0.18/-1.96/-3.47  &  0.18/-2.18/-3.90  &  \textcolor{red}{\textbf{0.20}}/\textcolor{red}{\textbf{-1.66}}/\textcolor{red}{\textbf{-2.13}}  &  0.19/-2.05/-3.80  \\
\hline
blonde+10\%  &  0.26/\textcolor{red}{\textbf{0.92}}/2.83  &  0.25/0.87/2.69  &  0.25/0.86/2.68  &  0.22/0.90/2.80  &  \textcolor{red}{\textbf{0.27}}/0.88/2.72  &  0.25/\textcolor{red}{\textbf{0.92}}/\textcolor{red}{\textbf{2.87}}  \\
blonde+30\%  &  0.25/0.82/2.57  &  0.26/0.76/2.21  &  0.26/0.70/2.14  &  0.25/0.80/2.47  &  \textcolor{red}{\textbf{0.27}}/0.76/2.35  &  \textcolor{red}{\textbf{0.27}}/\textcolor{red}{\textbf{0.85}}/\textcolor{red}{\textbf{2.59}}  \\
blonde+50\%  &  0.25/0.60/1.89  &  0.25/0.61/1.80  &  0.24/0.58/1.77  &  0.25/0.67/2.13  &  0.25/0.64/1.96  &  \textcolor{red}{\textbf{0.26}}/\textcolor{red}{\textbf{0.75}}/\textcolor{red}{\textbf{2.26}}  \\
blonde+70\%  &  0.24/0.02/0.22  &  0.22/0.35/1.02  &  0.21/0.35/1.05  &  0.22/0.33/1.11  &  0.20/0.38/1.19  &  \textcolor{red}{\textbf{0.26}}/\textcolor{red}{\textbf{0.49}}/\textcolor{red}{\textbf{1.47}}  \\
blonde+90\%  &  0.20/-0.81/-1.39  &  \textcolor{red}{\textbf{0.21}}/-0.62/-1.15  &  0.20/-0.65/-1.26  &  0.18/-0.75/-1.39  &  0.20/-0.52/-1.03  &  0.19/\textcolor{red}{\textbf{-0.45}}/\textcolor{red}{\textbf{-0.88}}  \\
\hline
\end{tabular}}
\end{center}
\end{table*}

\subsection{Impulse Noise Removal}\label{exp:inr}

\begin{figure*}[!t]
\captionsetup[subfigure]{justification=centering}
    \centering
     \begin{subfigure}{\imgwidl0tv}\fcolorbox{colorone}{colortwo}{\includegraphics[width=\textwidth,height=\imgheil0tv]{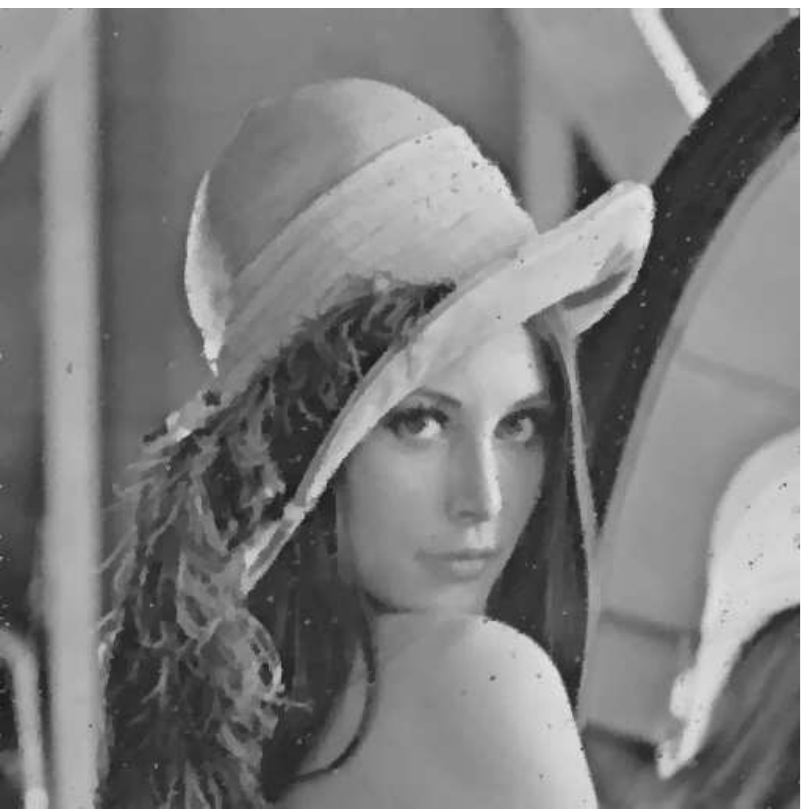}}\caption{\figlabsize by CVX\\~~~$SNR_0=0.92$}\end{subfigure}\ghs
      \begin{subfigure}{\imgwidl0tv}\fcolorbox{colorone}{colortwo}{\includegraphics[width=\textwidth,height=\imgheil0tv]{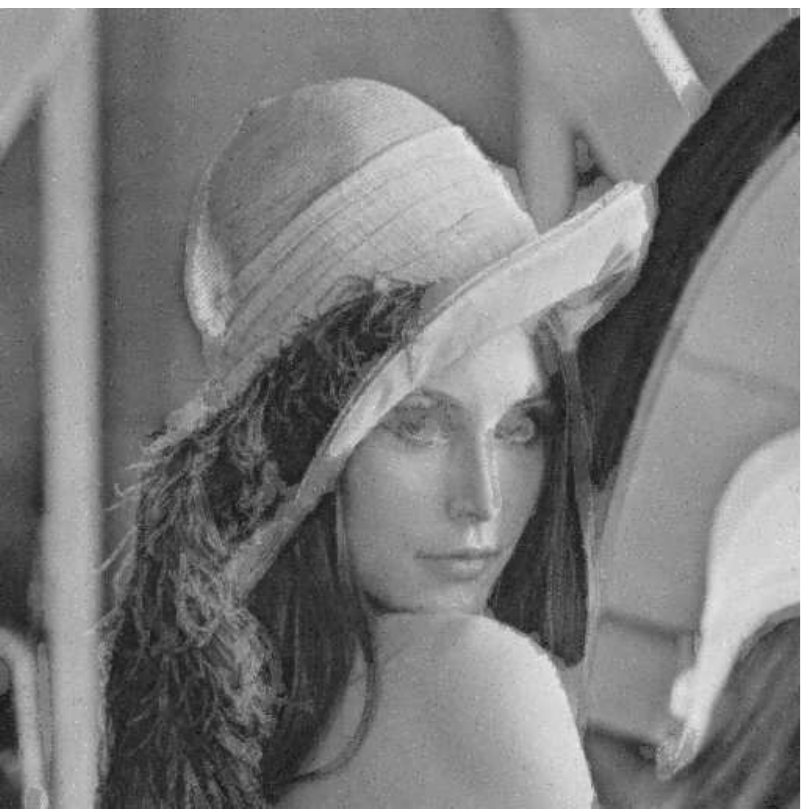}}\caption{\figlabsize by NI-ADM\\~~~$SNR_0=0.92$}\end{subfigure}\ghs
      \begin{subfigure}{\imgwidl0tv}\fcolorbox{colorone}{colortwo}{\includegraphics[width=\textwidth,height=\imgheil0tv]{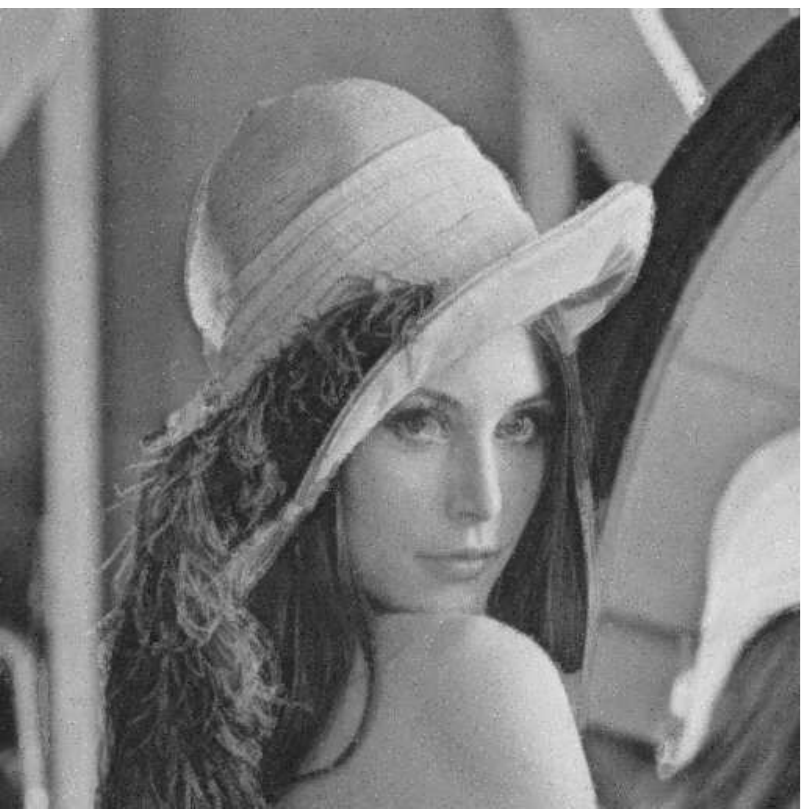}}\caption{\figlabsize by MD-ADM\\~~~$SNR_0=0.91$}\end{subfigure}

      \begin{subfigure}{\imgwidl0tv}\fcolorbox{colorone}{colortwo}{\includegraphics[width=\textwidth,height=\imgheil0tv]{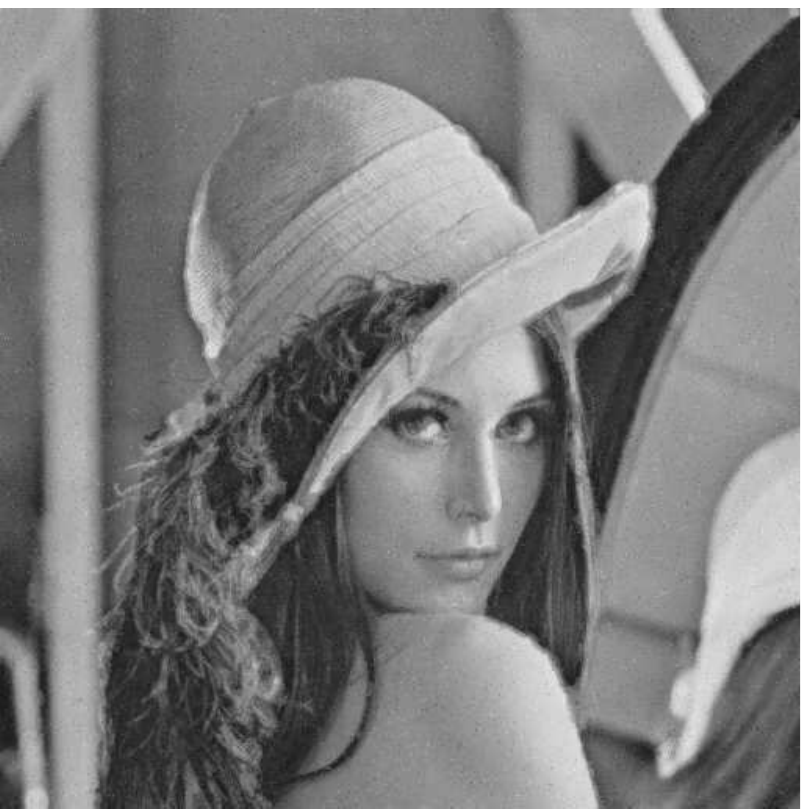}}\caption{\figlabsize by QPM\\~~~$SNR_0=0.94$}\end{subfigure}\ghs
      \begin{subfigure}{\imgwidl0tv}\fcolorbox{colorone}{colortwo}{\includegraphics[width=\textwidth,height=\imgheil0tv]{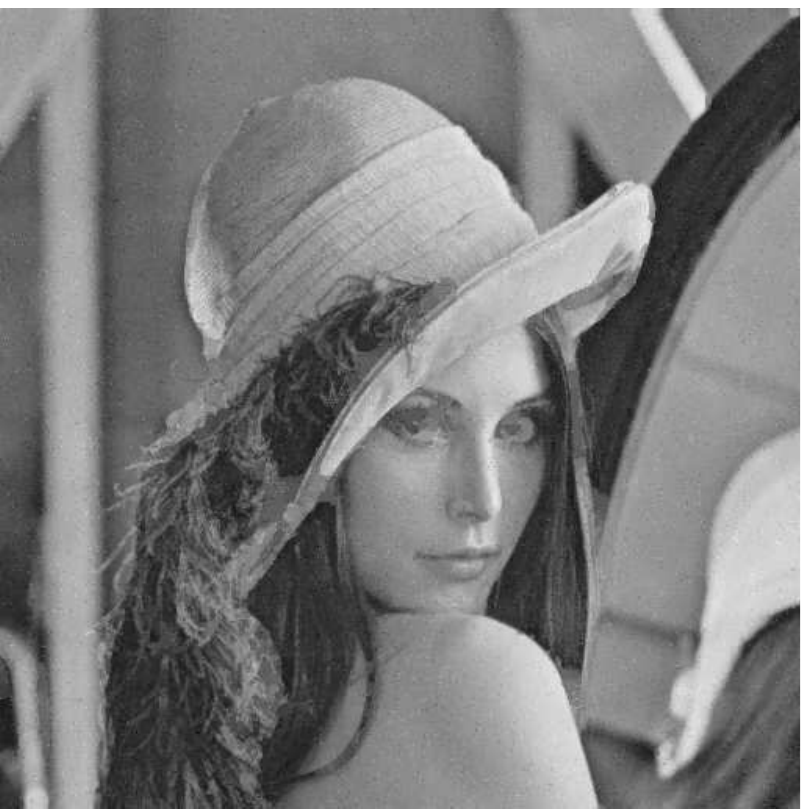}}\caption{\figlabsize by MPEC-EPM\\~~~$SNR_0=0.93$}\end{subfigure}\ghs
      \begin{subfigure}{\imgwidl0tv}\fcolorbox{colorone}{colortwo}{\includegraphics[width=\textwidth,height=\imgheil0tv]{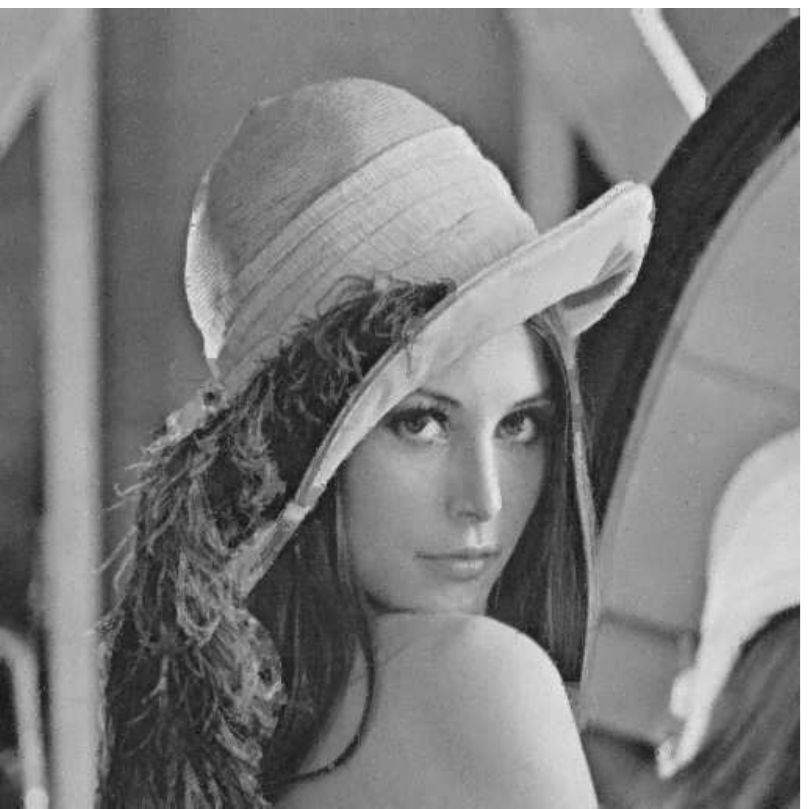}}\caption{\figlabsize by MPEC-ADM\\~~~$SNR_0=0.96$}\end{subfigure}
\vspace{0pt}
\caption{Random value impulse noise removal on `lenna' image.}
\label{fig:rv:removal}
\end{figure*}

Given a noisy image $\bbb{b}$ which was contaminated by random value impulsive noise, we consider a denoised image of $\bbb{x}$ as a minimizer of the following constrained $\ell_0TV$ model \cite{YuanG15}:
\beq
\min_{\bbb{x}}~TV(\bbb{x}),~s.t.~\|\bbb{x}-\bbb{b}\|_0\leq k
\eeq
\noi where $TV(\bbb{x}) \triangleq \|\bbb{\nabla} \bbb{x}\|_{p,1}$ denotes the total variation norm function, $\|\bbb{x}\|_{p,1}\triangleq\sum_{i=1}^n(|\bbb{x}_i|^p+|\bbb{x}_{n+i}|^p)^{\frac{1}{p}};~\bbb{\nabla} \triangleq \left[\bbb{\nabla}_x| \bbb{\nabla}_y\right]^T\in \mathbb{R}^{2n \times n}$, and $k$ specifies the noise level. Furthermore, $\bbb{\nabla}_x\in \mathbb{R}^{n\times n}$ and $\bbb{\nabla}_y\in \mathbb{R}^{n\times n}$ are two suitable linear transformation matrices that computes their discrete gradients. Following the experimental settings in \cite{YuanG15}, we use three ways to measure SNR, i.e. $SNR_0, SNR_1, SNR_2$. We show image recovery results in Table \ref{denoising:rv} when random-value impulse noise is added to the clean images \{`mandrill', `lenna', `lake', `jetplane', `blonde'\}. We also demonstrate a visualized example of random value impulse noise removal on `lenna' image in Figure \ref{fig:rv:removal}. A conclusion can be drawn that both MPEC-EPM and MPEC-ADM achieve state-of-the-art performance while MPEC-ADM is generally better than MPEC-EPM in this impulse noise removal problem.

\section{Conclusions}\label{sect:conc}
This paper presents two methods (exact penalty and alternating direction) to solve the general sparsity constrained minimization problem. Although it is non-convex, we design effective penalization/regularization schemes to solve its equivalent MPEC. We also prove that both of our methods are convergent to first-order KKT points. Experimental results on a variety of sparse optimization applications demonstrate that our methods achieve state-of-the-art performance.





\vspace{15pt}
\appendix

\noi {\LARGE \textbf{Appendix}}
\vspace{-2pt}

\section{Matlab Code of Breakpoint Search Algorithm\label{supp:breakpoint}}

In this section, we include a Matlab implementation of the breakpoint search algorithm \cite{HelgasonKL80} which solves the optimization problem in Eq (\ref{eq:clap_simplex}) exactly in $n\log(n)$ time.
\beq \label{eq:clap_simplex}
\min_{\bbb{x}}~\frac{1}{2}\bbb{x}^Tdiag(d)\bbb{x} + \bbb{x}^T\bbb{a},~s.t.~\bbb{0}\leq \bbb{x}\leq\bbb{1},~\bbb{x}^T\bbb{1} ~\lesseqqgtr ~s
\eeq
\noi where $\bbb{x},\bbb{d},\bbb{a}\in \mathbb{R}^n$, $\bbb{d}>0$, and $0\leq s\leq n$. In addition, $\lesseqqgtr$ is a comparison operator, and it can be strict equal operator (=), greater than or equal operator ($\geq$), or less than or equal operator ($\leq$). We list our code as follows.


\begin{lstlisting}
function [x] = BreakPointSearch(d,a,s,oper)
% This program solves the following QP: 
% min_x 0.5*x'*diag(d)*x + a'*x
% s.t. sum(x) oper s,  0 <= x <= 1
% oper can be '==', '>=', or '<='
% We assume d > 0.

n = length(d);assert( s>=0 && s<=n );
y = sort([-a-d;-a],'ascend');
l = 1; r = 2*n; L = n; R = 0;
while 1
    if(r-l==1)
        S = (s-L)/(R-L);
        rho = y(l)+(y(r) - y(l))*S;
        break;
    else
        m = floor(0.5*(l+r));
        t = (-a-y(m))./d;
        C = sum(max(0,min(t,1)));
        if(C==s),
            rho = y(m);break;
        elseif(C>s)
            l = m;L=C;
        elseif(C<s)
            r = m;R=C;
        end
    end
end
if (oper(1)=='>')
rho = min(rho,0); 
elseif(oper(1)=='<')
rho = max(rho,0); 
end
t = (-a-rho)./d; 
x = max(0,min(t,1));

\end{lstlisting}

%

\fontsize{11.5}{12}\selectfont
\bibliography{mybib}
\end{document}